\documentclass[12pt, a4paper]{article}
\title{A note on the bicategory of Landau-Ginzburg models ($\mathcal{LG}_K$)}
\author{ Yves Baudelaire Fomatati$^{a}$\footnote{a. Department of Mathematics and Statistics, University of Ottawa, yfomatat@uottawa.ca .
 }}

\date{}

\usepackage{xcolor}

\usepackage[square,comma,numbers,sort&compress]{natbib}
\usepackage[centertags]{amsmath}
\usepackage{latexsym}
\usepackage{leqno}
\usepackage{amsfonts}
\usepackage{amsmath,amssymb}
\usepackage{enumerate}

\usepackage{amsthm}
\usepackage[utf8,latin9]{inputenc}

\theoremstyle{plain}
\newtheorem{remark}{Remark}[section]
\theoremstyle{plain}
\newtheorem{lemma}{Lemma}[section]
\theoremstyle{plain}
\newtheorem{proposition}{Proposition}[section]
\theoremstyle{plain}
\newtheorem{theorem}{Theorem}[section]
\theoremstyle{plain}
\newtheorem{definition}{Definition}[section]
\theoremstyle{plain}

\theoremstyle{plain}

\theoremstyle{plain}
\newtheorem{notation}{Notations}[section]
\newtheorem{example}{Example}[section]

\theoremstyle{plain}

\usepackage[all, 2cell]{xy}
\usepackage{txfonts}

\begin{document}
\maketitle
\begin{quote}
  \textbf{Abstract}
\end{quote}
The bicategory of Landau-Ginzburg models denoted by $\mathcal{LG}_{K}$ possesses adjoints and this helps in explaining a certain duality that exists in the setting of Landau-Ginzburg models in terms of some specified relations. The construction of $\mathcal{LG}_{K}$ is reminiscent of, but more complex than, the construction of the bicategory of associative algebras and bimodules. In this paper,
we review this complex but very inspiring construction in order to expose it more to pure mathematicians. 
In particular, we spend some time explaining the intricate construction of unit morphisms in this bicategory from a new vantage point. Besides, we briefly discuss how this bicategory could be constructed in more than one way using the variants of the Yoshino tensor product. Furthermore, without resorting to Atiyah classes, we prove that the left and right unitors in this bicategory have direct right inverses but do not have direct left inverses.

\section{Introduction}
A Landau-Ginzburg model is a model in solid state physics for superconductivity. Superconductivity is a quantum mechanical phenomenon of exactly zero electrical resistance and expulsion of magnetic flux field occurring in superconductors when cooled below their critical temperature (i.e., a temperature at and above which their vapor cannot be liquefied no matter how much pressure is applied). Some excellent textbooks  and papers giving a detailed account on superconductivity are \cite{tilley2019superfluidity} and \cite{onnes1938akad}.\\
Landau-Ginzburg models play a significant role in several areas of pure
mathematics and mathematical physics including singularity theory, representation theory, and conformal or topological field theory. The connection between these areas motivates the study of Landau-Ginzburg models.
In \cite{carqueville2016adjunctions}, it is explained how Landau-Ginzburg models give rise to a bicategory with adjoints (also called duals) and the structure maps in this bicategory are described, which include the units and counits of adjunction in terms of basic invariants called \textit{Atiyah classes} (section 3 of \cite{carqueville2016adjunctions}). The authors of \cite{carqueville2016adjunctions} studied this bicategory of Landau-Ginzburg models and showed that the bicategorical perspective offers a certain unified approach to Landau-Ginzburg models. They used Atiyah classes to prove that the unitors in this bicategory have homotopy inverses . \\
 In this paper, we elucidate the intricate construction of the unit morphisms and without resorting to Atiyah classes, we prove that the unitors have direct right inverses. Unfortunately, they do not have direct left inverses.
As described in \cite{carqueville2016adjunctions}, this bicategory has \textit{potentials} (definition 2.4 of \cite{carqueville2016adjunctions}) as objects and matrix factorizations as $1-$morphisms.
\\
But, unlike \cite{carqueville2016adjunctions} we do not restrict ourselves to \textit{potentials} which are polynomials satisfying some conditions (definition 2.4 p.8 of \cite{carqueville2016adjunctions} )\footnote{For an earlier reference, see section 3 p.17-18 of \cite{khovanov2008matrix}.}. Potentials have many applications in quantum algebra \cite{khovanov2008matrix} and mathematical physics \cite{carqueville2016adjunctions}.
It turns out that the authors of \cite{carqueville2016adjunctions} used potentials to suit their purposes because even if we take the objects of $\mathcal{LG}_{K}$ to be polynomials rather than potentials and then apply the construction of $\mathcal{LG}_{K}$ given in \cite{carqueville2016adjunctions}, we obtain virtually the same bicategory except that we now have more objects. Moreover, in this new category whose objects are polynomials (without restrictions), the unit construction of the bicategory of Landau-Ginzburg models is done the same way. Therefore, in our presentation, the objects of $\mathcal{LG}_{K}$ are polynomials in general.\\
 The composition of $1-$morphisms in $\mathcal{LG}_{K}$ is done using the Yoshino tensor product of matrix factorizations. We are going to observe that in place of this product, any of its variants could also be used thereby yielding variants of the bicategory $\mathcal{LG}_{K}$.
\\

In the next section, we give a short note on linear factorizations and in section 3, the construction of $\mathcal{LG}_K$ is reviewed.
\\Except otherwise stated, $R=K[x]$ where $x=x_{1},\cdots, x_{n}$ and $K$ is a commutative ring with unity.
\section{Matrix factorizations and the Yoshino tensor product}
In this section, we set the stage for the review of the bicategory $\mathcal{LG}_K$ (cf. section \ref{bicategory LG}). First, we give some useful definitions and make some observations. Next, we recall Yoshino tensor product of matrix factorizations and its variants. We also recall the notion of homotopic linear factorizations. In the third subsection, we briefly discuss factorizations of finite rank. Finally, in the last subsection we revisit the tensor product of matrix factorizations. In fact, after recalling two different ways in which a matrix factorization can be defined, we show how the Yoshino tensor product can be used on matrix factorizations of two polynomials to produce a matrix factorization of their sum.
\subsection{Some definitions and observations}\label{link lin facto n complexes}

\begin{definition}(p.8 of \cite{carqueville2016adjunctions}) Linear factorization \label{def lin facto}\\
A \textbf{linear factorization} of $f\in R$ is a $\mathbb{Z}_{2}-$graded $R-$module $X=X^{0}\oplus X^{1}$
together with an odd (i.e., grade reversing) $R-$linear endomorphism $d: X\longrightarrow X$
such that $d^{2}=f\cdot id_{X}$.\\
$f\cdot id_{X}$ stands for the endomorphism $x\mapsto f\cdot x,\,\,\forall x\in X$.
\end{definition}
$d$ is called a \textit{twisted differential} in \cite{carqueville2015toolkit}.
 $d$ is actually a pair of maps $(d^{0},d^{1})$ that we may depict as follows:
$$\xymatrix {X^{0}\ar [r]^{d^{0}} &X^{1}\ar [r]^{d^{1}} & X^{0}}$$
such that $d^{0}\circ d^{1}=f\cdot id_{X^{1}}$ and $d^{1}\circ d^{0}= f\cdot id_{X^{0}}$.

\begin{example}\label{exple of lin facto}
Keeping to the notations of definition \ref{def lin facto}, let $R=\mathbb{C}[x]$ consider $f=x^{3}\in R$.
Then take $X$ to be the $\mathbb{Z}_{2}$-graded $\mathbb{C}[x]$-module
$\mathbb{C}[x]\oplus \mathbb{C}[x]$. Then define $d: X\rightarrow X$ as follows:

$$\xymatrix {X^{0}\ar [r]^{M} &X^{1}\ar [r]^{M} & X^{0}}$$

where some of the choices for $M$ are the following:

\[
M=
  \begin{bmatrix}
    0 & x  \\
    x^{2} & 0
  \end{bmatrix}
\]
or

\[
M=
  \begin{bmatrix}
    0 & 1  \\
    x^{3} & 0
  \end{bmatrix}
\]
Here, $d^{0}=d^{1}$ and $M$ is the matrix corresponding to the $R$-linear endomorphism $d^{0}$.
In general, if $f=x^{n}$, then we can take $M$ to be
\[
M_{q}=
  \begin{bmatrix}
    0 & x^{q}  \\
    x^{n-q} & 0
  \end{bmatrix}
\]
Clearly $d^{2}=f\cdot id_{X}$.
\end{example}

\begin{remark}
  If $X$ is a free $R-$module, then the pair $(X,d)$ is called a \textit{matrix factorization}. If $M_{0}$ and $M_{1}$ are respectively matrices of the $R-$linear endomorphisms $d^{0}$ and $d^{1}$, then the pair $(M_{0},M_{1})$ would be a matrix factorization of $f$ according to definition 2.1 of \cite{fomatati2019multiplicative}.
\end{remark}
Eisenbud \cite{eisenbud1980homological} was the first to introduce the notion of matrix factorizations. He related them to maximal Cohen-Macaulay modules
(cf. ( Chap. 7 of \cite{yoshino1990maximal}), \cite{schreyer1987finite}).\\
The original definition of a matrix factorization given by Eisenbud (p.15 of \cite{eisenbud1980homological}) is as follows: a matrix factorization of an element $f$ in a ring $R$ (with unity) is an ordered pair of maps of free $R-$modules $\phi: F\rightarrow G$ and $\psi: G \rightarrow F$ s.t., $\phi\psi=f\cdot 1_{G}$ and $\psi\phi=f\cdot 1_{F}$. From this definition of Eisenbud, it is clear that to obtain a linear factorization from a matrix factorization, we need to have $F$ and $G$ represent respectively the even and the odd components of a $\mathbb{Z}_{2}-$graded $R-$module $X$ such that $\phi$ and $\psi$ are grade reversing maps. Just like in \cite{carqueville2016adjunctions}, in this section, we will often refer to a matrix factorization $(X,d_{X})$ by $X$ without mentioning the differential $d_{X}$.

Given a basis for $X$, it is sometimes more convenient to identify $d_{X}$ with its associated block matrix
$$ d_{X}= \begin{pmatrix}
0 & d_{X}^{1} \\
d_{X}^{0} & 0
\end{pmatrix}$$

Two remarkable operations can be carried out on matrix factorizations, namely \textit{direct sum} and \textit{tensor product}.
The direct sum of two matrix factorizations $X$ and $Y$ is defined in the obvious way:
$$(X\oplus Y)^{i}= X^{i} \oplus Y^{i}\;\; and \;\; d^{i}_{X\oplus Y}=d^{i}_{X} + d^{i}_{Y}$$
\\ As for the operation of tensor product it is presented in the following subsection.

\subsection{Yoshino's tensor product of matrix factorization and its variants} \label{subsec: Yoshino tens prod and variant}
Here, we recall the Yoshino tensor product of matrix factorization and its variants. First, the definition of the category of matrix factorization of a power series is recalled.
\begin{definition} \cite{yoshino1998tensor}  \label{defn matrix facto of polyn}   \\
An $n\times n$ \textbf{matrix factorization} of a power series $f\in \;R$ is a pair of $n$ $\times$ $n$ matrices $(\phi,\psi)$ such that
$\phi\psi=\psi\phi=fI_{n}$, where $I_{n}$ is the $n \times n$ identity matrix and the coefficients of $\phi$ and of $\psi$ are taken from $R$.
\end{definition}
Also recall ($\S 1$ of \cite{yoshino1998tensor}) the definition of
the category of matrix factorizations of a power series $f\in R=K[[x]]:=K[[x_{1},\cdots,x_{n}]]$ denoted by $MF(R,f)$ or $MF_{R}(f)$, (or even $MF(f)$ when there is no risk of confusion):\\
$\bullet$ The objects are the matrix factorizations of $f$.\\
$\bullet$ Given two matrix factorizations of $f$; $(\phi_{1},\psi_{1})$ and $(\phi_{2},\psi_{2})$ respectively of sizes $n_{1}$ and $n_{2}$, a morphism from $(\phi_{1},\psi_{1})$ to $(\phi_{2},\psi_{2})$ is a pair of matrices $(\alpha,\beta)$ each of size $n_{2}\times n_{1}$ which makes the following diagram commute \cite{yoshino1998tensor}:
$$\xymatrix@ R=0.6in @ C=.75in{K[[x]]^{n_{1}} \ar[r]^{\psi_{1}} \ar[d]_{\alpha} &
K[[x]]^{n_{1}} \ar[d]^{\beta} \ar[r]^{\phi_{1}} & K[[x]]^{n_{1}}\ar[d]^{\alpha}\\
K[[x]]^{n_{2}} \ar[r]^{\psi_{2}} & K[[x]]^{n_{2}}\ar[r]^{\phi_{2}} & K[[x]]^{n_{2}}}$$
That is,
$$\begin{cases}
 \alpha\phi_{1}=\phi_{2}\beta  \\
 \psi_{2}\alpha= \beta\psi_{1}
\end{cases}$$
More details on this category are found in chapter 2 of \cite{fomatati2019multiplicative}.
\\
\begin{definition} \cite{yoshino1998tensor} \label{defn Yoshino tensor prodt}
Let $X=(\phi,\psi)$ be an $n\times n$ matrix factorization of $f\in R$  and $X'=(\phi',\psi')$ an $m\times m$ matrix factorization of $g\in S$. These matrices can be considered as matrices over $L=K[[x,y]]$ and the \textbf{tensor product} $X\widehat{\otimes} X'$ is given by\\
(\(
\begin{bmatrix}
    \phi\otimes 1_{m}  &  1_{n}\otimes \phi'      \\
   -1_{n}\otimes \psi'  &  \psi\otimes 1_{m}
\end{bmatrix}
,
\begin{bmatrix}
    \psi\otimes 1_{m}  &  -1_{n}\otimes \phi'      \\
    1_{n}\otimes \psi'  &  \phi\otimes 1_{m}
\end{bmatrix}
\))\\
where each component is an endomorphism on $L^{n}\otimes L^{m}$.
\end{definition}
 $X\widehat{\otimes} X'$ is in fact an object of $MF_{L}(f+g)$ of size $2nm$.

\begin{remark}
When $n=1$, we get a $1$ $\times$ $1$ matrix factorization of $f$, i.e., $f=[g][h]$ which is simply a factorization of $f$ in the classical sense. But in case $f$ is not reducible, this is not interesting, that is why we will mostly consider $n > 1$.
\end{remark}
\textbf{Variants of Yoshino's tensor product of matrix factorizations}\\
\begin{definition} \label{defn Yoshino tensor prodt} \cite{fomatati2022tensor}
Let $X=(\phi,\psi)$ be an $n\times n$ matrix factorization of $f\in R$  and $X'=(\phi',\psi')$ an $m\times m$ matrix factorization of $g\in S$. These matrices can be considered as matrices over $L=K[[x,y]]$ and the \textbf{tensor products} $X\widehat{\otimes}' X'$, $X\widehat{\otimes}'' X'$ and $X\widehat{\otimes}''' X'$ are respectively given by\\
(\(
\begin{bmatrix}
  1_{n}\otimes \phi'  &  \psi\otimes 1_{m}  \\
    \phi\otimes 1_{m} & -1_{n}\otimes \psi'
\end{bmatrix}
,
\begin{bmatrix}
   1_{n}\otimes \psi' & \psi\otimes 1_{m}         \\
    \phi\otimes 1_{m} &  -1_{n}\otimes \phi'
\end{bmatrix}
\)),\\
(\(
\begin{bmatrix}
  \psi\otimes 1_{m} & -1_{n}\otimes \psi' \\
  1_{n}\otimes \phi' & \phi\otimes 1_{m}
\end{bmatrix}
,
\begin{bmatrix}
 \phi\otimes 1_{m} & 1_{n}\otimes \psi' \\
 -1_{n}\otimes \phi' & \psi\otimes 1_{m}
\end{bmatrix}
\)) and \\
(\(
\begin{bmatrix}
-1_{n}\otimes \psi' & \phi\otimes 1_{m}\\
\psi\otimes 1_{m}  & 1_{n}\otimes \phi'
\end{bmatrix}
,
\begin{bmatrix}
-1_{n}\otimes \phi' & \phi\otimes 1_{m} \\
 \psi\otimes 1_{m} & 1_{n}\otimes \psi'
\end{bmatrix}
\))\\
where each component is an endomorphism on $L^{n}\otimes_{L} L^{m}$.
\end{definition}
 $X\widehat{\otimes}' X'$, $X\widehat{\otimes}'' X'$ and $X\widehat{\otimes}''' X'$ are in fact objects of $MF_{L}(f+g)$, each of size $2nm$. These objects are mutually distinct.
\begin{definition} \cite{fomatati2022tensor}
  $\widehat{\otimes}'$, $\widehat{\otimes}''$ and $\widehat{\otimes}'''$ are respectively called the first, second and third variant of the Yoshino tensor product.
\end{definition}

\begin{proposition} \cite{fomatati2022tensor}
  $\widehat{\otimes}'$, $\widehat{\otimes}''$ and $\widehat{\otimes}'''$ are functorial operations.
\end{proposition}




\begin{definition} (p.9 \cite{carqueville2016adjunctions}) Morphism of linear factorizations \label{morphism of matrix facto}\\
  A \textbf{morphism of linear factorizations} $(X,d_{X})$ and $(Y,d_{Y})$ is an even (i.e., a grade preserving) $R-$linear map $\phi: X \longrightarrow Y$ such that $d_{Y}\phi = \phi d_{X}$.
\end{definition}
Concretely (see page 19 of \cite{khovanov2008matrix}), $\phi$ is a pair of maps $\xymatrix {X^{0}\ar [r]^{\phi^{0}} &Y^{0}}$ and $\xymatrix {X^{1}\ar [r]^{\phi^{1}} &Y^{1}}$ such that the following diagram commutes:\\

 $$\xymatrix@ R=0.6in @ C=.75in{X^{0} \ar[r]^{{d^{0}_{X}}} \ar[d]_{\phi^{0}} &
X^{1} \ar[d]^{\phi^{1}} \ar[r]^{{d^{1}_{X}}} & X^{0}\ar[d]^{\phi^{0}}\\
Y^{0} \ar[r]^{{d^{0}_{Y}}} & Y^{1}\ar[r]^{{d^{1}_{Y}}} & Y^{0}}$$

\begin{remark}
 The family of $R-$linear maps $Hom_{R}(X,Y)$ between two linear factorizations $X$ and $Y$ is an $R-$module with action $r(\phi^{0}, \phi^{1})= (r\phi^{0},r\phi^{1})$, $r\in R$.\\
Since $X$ and $Y$ are linear factorizations of $f\in R$, they are $\mathbb{Z}_{2}-$graded modules, thus we can write $X=X^{0}\oplus X^{1}$ and $Y=Y^{0}\oplus Y^{1}$. Hence maps in $Hom_{R}(X,Y)$ are of degree one or zero. Thus, we can write $Hom_{R}(X,Y)= Hom^{1}_{R}(X,Y)\oplus Hom^{0}_{R}(X,Y)$.
\end{remark}


\begin{proposition} \label{Hom(X,Y) is a Z2-graded complex}
Let $X$ and $Y$ be two linear factorizations of an element $f\in R$. Then $Hom_{R}(X,Y)$ is a $\mathbb{Z}_{2}-$graded complex with differential $$d: \phi\mapsto d_{Y}\circ \phi -(-1)^{|\phi|} \phi \circ d_{X} $$
where $|\phi|$ is the degree
of the map $\phi$.
\end{proposition}
N.B. we will drop the $R$ in $Hom_{R}(X,Y)$ for ease of notation.
\begin{proof} see Prop. 5.1 of \cite{fomatati2019multiplicative}.

\end{proof}
One of the main goals of this paper is to review the construction of the bicategory $\mathcal{LG}_{K}$ and in particular to elucidate the intricate construction of the unit in this bicategory. To that end, we need more ingredients. Besides the notion of tensor products of matrix factorizations we will need the notion of homotopy between linear factorizations.
\begin{notation} \label{notation F(R,f) is catg of lin facto}
Let $R$ be a commutative ring and $f\in R$. We will write $F(R,f$) for the category whose objects are
linear factorizations of $f$ and morphisms are homomorphisms of linear factorizations.
\end{notation}

\begin{definition} \label{homotopic lin facto} (p.9 \cite{carqueville2016adjunctions}) homotopic linear factorizations\\
 Let $(X,d_{X})$ and $(Y,d_{Y})$ be linear factorizations. Two morphisms $\varphi, \psi : X\longrightarrow Y$ are \textbf{homotopic} if there exists an odd $R-$linear map $\lambda: X\longrightarrow Y$ such that $d_{Y}\lambda + \lambda d_{X}= \psi - \varphi$.\\
 More precisely, the following diagram commutes:
$$\xymatrix@ R=0.6in @ C=.75in{X_{1} \ar[r]^{d_{X}} \ar[d]_{\psi-\phi} &
X_{0}\ar[dl]_{\lambda_{0}} \ar[d]^{\psi-\phi} \ar[r]^{d_{X}} & X_{1}\ar[dl]_{\lambda_{1}}\ar[d]^{\psi-\phi \,\,\,\,\,\,\,\,\,\,\,\,\,\,\,\dag}\\
Y_{1} \ar[r]^{d_{Y}} & Y_{0}\ar[r]^{d_{Y}} & Y_{1}}$$ i.e., $$d_{Y}\circ \lambda_{0} + \lambda_{1}\circ d_{X}=\psi-\phi$$
\end{definition}

Recall that given a category $C$, a \textit{congruence relation} $R$ on $C$, is an equivalence relation $R_{XY}$ on $Hom_{C}(X,Y)$ for objects $X,Y$ s.t. the equivalence relations are compatible with composition. 
Moreover, we know
that given a congruence relation $\mathcal{R}$ on $C$, we can define the quotient category $C/\mathcal{R}$ whose objects are those of $C$ but morphisms are equivalence classes of morphisms in $C$. Composition of morphisms in $C/\mathcal{R}$ is well defined since $\mathcal{R}$ is a congruence relation.\\
It is easy to see that equality up to homotopy is an equivalence relation. This relation is compatible with composition and so one can form the quotient category denoted by $HF(R,f$) in \cite{carqueville2016adjunctions}.

When defining a matrix factorization in the previous section, we did not talk about the rank of a matrix factorization, it is time to do so.

\subsection{A glimpse of factorizations of finite rank}
This subsection mostly relies on work done in \cite{khovanov2008matrix} (p.18, 20) and \cite{carqueville2016adjunctions}). It is reproduced here for the sake of completeness.\\
If $X$ is a matrix factorization, then it is not always the case that the ranks of the $R-$modules $X^{0}$ and $X^{1}$ be finite. This notion of finiteness of the rank of a matrix factorization was crucial in the study of Morita contexts in the bicategory $\mathcal{LG}_{K}$ \cite{fomatati2021necessary}.

\begin{definition}\label{finite rank matrix facto}
We say that a matrix factorization is of \textit{finite rank} if its underlying free $R$-module is of \textbf{finite rank}.
\end{definition}

Given a basis for the $R$-module $X=X^{0}\oplus X^{1}$, the differential $d_{X}$ is sometimes identified with the matrix $$d_{X}= \begin{pmatrix}
0 & d_{X}^{1} \\
d_{X}^{0} & 0
\end{pmatrix}$$
\\

 Thus, $$d^{2}_{X}= \begin{pmatrix}
0 & d_{X}^{1} \\
d_{X}^{0} & 0
\end{pmatrix}\begin{pmatrix}
0 & d_{X}^{1} \\
d_{X}^{0} & 0
\end{pmatrix}
=
\begin{pmatrix}
d_{X}^{1}d_{X}^{0}& 0\\
0 & d_{X}^{0}d_{X}^{1}
\end{pmatrix}
=\begin{pmatrix}
f\cdot id_{X^{0}} & 0 \\
0 & f\cdot id_{X^{1}}
\end{pmatrix}
=
f\begin{pmatrix}
 id_{X^{0}} & 0 \\
0 &  id_{X^{1}}
\end{pmatrix}
$$
So, $d^{2}_{X}=fId$ where $Id= \begin{pmatrix}
 id_{X^{0}} & 0 \\
0 &  id_{X^{1}}
\end{pmatrix}$.

\begin{notation}\label{notations for categ of facto} We keep the following notations used in \cite{carqueville2016adjunctions}:\\
$\bullet$ As alluded above, $HF(R, f)$ denotes the category of linear factorizations of $f\in R$ modulo homotopy. \\
$\bullet$ We also denote by $HMF(R, f)$ its full subcategory of matrix factorizations. \\
$\bullet$ Furthermore, we write $hmf(R, f)$ for the full subcategory of finite rank matrix factorizations, viz. the matrix factorizations whose underlying $R$-module is free of finite rank.\\
Recall that the category of matrix factorizations of $f$ is denoted $MF(R, f)$ (or sometimes simply denoted by $MF(f)$ when there is no risk of confusion).
\end{notation}
  $$Hom_{HMF(R, f)}(X,Y)= Hom_{MF(R, f)}(X,Y)/\{Null-homotopic\;maps\}$$
A null-homotopic map is a map that is homotopic to the zero map.
$MF(R, f)$ is additive and $R-$linear (p.19 \cite{khovanov2008matrix})\footnote{In \cite{khovanov2008matrix}, the author writes $MF^{all}_{f}$ instead of $MF(R, f)$}.

If we choose bases in the free $R-$modules $X^{0}$ and $X^{1}$, then we can write the maps $d^{0}_{X}$ and $d^{1}_{X}$ as $m\times m$ matrices $D^{0}_{X}$ and $D^{1}_{X}$ with coefficients in $R$. These matrices satisfy the equalities: $$D^{0}_{X}D^{1}_{X}=f\cdot Id,\;\;\;D^{1}_{X}D^{0}_{X}=f\cdot Id$$ where $Id$ stands for the identity $m\times m$ matrix.
Alternatively, we can describe this factorization by a $2m\times 2m$ matrix with off-diagonal blocks $D^{0}$ and $D^{1}$:
$$D= \begin{pmatrix}
0 & D^{1} \\
D^{0} & 0
\end{pmatrix}\;\;\;\;\;D^{2}=f\cdot Id,\,or\,written\,simply\,D^{2}=fId,$$
where we dropped the subscripts on $D$ for ease of notation. The $Id$ in the last equality above is the identity matrix of size $2m$.\\
Matrix description of objects in $hmf(R,f)$ extends to infinite rank factorizations.
Matrices $D^{0}$ and $D^{1}$ then have infinite rank, but each of their columns has only finitely many non-zero entries.
If factorizations $X$ and $Y$ are written in matrix form, $X=(D^{0}_{X},D^{1}_{X})$ and $Y=(D^{0}_{Y},D^{1}_{Y})$, then a homomorphism $g: X\longrightarrow Y$ is a pair of matrices $(G_{0}, G_{1})$
such that $G_{1}D^{0}_{X}= D^{0}_{Y}G_{0}$ and $G_{0}D^{1}_{X}=D^{1}_{Y}G_{1}$.

\begin{remark}\label{eqce of equations of facto}
  The two equations $G_{1}D^{0}_{X}= D^{0}_{Y}G_{0}$ and $G_{0}D^{1}_{X}=D^{1}_{Y}G_{1}$ are equivalent (p.20 \cite{khovanov2008matrix}, p.174 \cite{leuschke2012cohen}). We give a proof of the equivalence here.
\\
\begin{eqnarray*}
G_{0}D^{1}_{X}=D^{1}_{Y}G_{1}&\Rightarrow & G_{0}D^{1}_{X}D^{0}_{X}=D^{1}_{Y}G_{1}D^{0}_{X}\\
&\Rightarrow& D^{0}_{Y}G_{0}D^{1}_{X}D^{0}_{X}=D^{0}_{Y}D^{1}_{Y}G_{1}D^{0}_{X}\\
&\Rightarrow& D^{0}_{Y}G_{0}f \cdot Id = f\cdot Id G_{1}D^{0}_{X}\\
&\Rightarrow& f \cdot D^{0}_{Y}G_{0} Id = f\cdot Id G_{1}D^{0}_{X}\\
&\Rightarrow& f \cdot Id D^{0}_{Y}G_{0} = f\cdot Id G_{1}D^{0}_{X}\;\;since\;D^{0}_{Y}G_{0} Id=Id D^{0}_{Y}G_{0}\\
&\Rightarrow& (f \cdot Id) (D^{0}_{Y}G_{0} - G_{1}D^{0}_{X})=0\\
&\Rightarrow& D^{0}_{Y}G_{0}=G_{1}D^{0}_{X} \;\;\;\;\;\;\;\; as \;\;desired.
\end{eqnarray*}

\end{remark}
The third implication above holds because $D^{1}_{X}D^{0}_{X}=f\cdot Id = D^{0}_{Y}D^{1}_{Y}$. The fourth holds because $f\in R$ and so can be moved across matrices. The last implication is true because here $f$ is an arbitrary element of $R$ and $f\cdot Id$ cannot be zero for all $f$.
Similarly, one can show that $D^{0}_{Y}G_{0}=G_{1}D^{0}_{X}\Rightarrow G_{0}D^{1}_{X}=D^{1}_{Y}G_{1}$.
Hence the equivalence holds.

\subsection{Tensor Products of Matrix Factorizations (Revisited)} \label{subsec:tensor prod}
In this subsection, we show how to use the Yoshino tensor product of matrix factorizations to produce a factorization of sum of two polynomials when their respective factorizations are given in the form $(X,d_{X})$ as defined in Definition \ref{def lin facto}. In addition, we show that variants of the Yoshino tensor product can be used to produce new factorizations which all differ at the level of their differentials. We explain how to obtain these differentials. This explanation helps understand how the composition of $1-$morphisms in the bicategory $\mathcal{LG}_{K}$ could be performed using variants of the Yoshino tensor product. As discussed in section \ref{bicategory LG}, this gives rise to variants of $\mathcal{LG}_{K}$.

The notion of tensor product of matrix factorizations first appeared in Yoshino's paper \cite{yoshino1998tensor}. Recall that Yoshino \cite{yoshino1998tensor} defines an $n\times n$ matrix factorization of a polynomial $f$ to be a pair of matrices $(P,Q)$ such that $PQ=f\cdot I_{n}$ (see definition \ref{defn matrix facto of polyn}).
\\N.B. Though Yoshino's tensor product is written with a hat, henceforth we will sometimes drop the hat of $\widehat{\otimes}$ and simply write $\otimes$ for ease of notation whenever there will be no risk of confusion. Likewise, we will also sometimes write its three variants without their hats $\otimes '$, $\otimes ''$ and $\otimes '''$. \\
Section 2 of \cite{carqueville2016adjunctions} gives another definition of matrix factorizations in terms of $\mathbb{Z}_{2}$-graded free modules.
In fact, a matrix factorization of $f\in R=K[x_{1},\cdots,x_{n}]$ can also be defined as a $\mathbb{Z}_{2}-$graded free $R-$module $X=X^{0}\oplus X^{1}$ together with an odd $R-$linear endomorphism also called differential map, $d_{X}$ s.t. $d_{X}^{2}=f\cdot 1_{X}$ (cf. section 2 of \cite{carqueville2016adjunctions}). Thus, if we consider the pair of matrices representing the (even and odd) components of such a differential map, we recover the definition of a matrix factorization of $f$ as earlier recalled (cf. definition \ref{defn matrix facto of polyn}).\\
 Henceforth in this paper, we will often simply write "factorization(s)" in place of "matrix factorization(s)".\\
  We are soon going to recall (cf. section \ref{bicategory LG}) that a $1-$morphism in the bicategory $\mathcal{LG}_{K}$ is a factorization of the difference of two polynomials. It is important to recall that the definition of matrix factorization used in \cite{carqueville2016adjunctions} is the one given in section 2 of \cite{carqueville2016adjunctions}, i.e., a pair $(X,d_{X})$ where $X$ is a graded module and $d_{X}$ is the differential map as explained above. We saw in definition \ref{defn Yoshino tensor prodt} how to construct a factorization of the sum of two polynomials from their respective factorizations. But also remark that in definition \ref{defn Yoshino tensor prodt}, the factorizations are given as pairs of matrices and not in the form $(X,d_{X})$ as defined in \cite{carqueville2016adjunctions}.\\
   So, the question is: how do we use the Yoshino tensor product on say $(X,d_{X})$ (a factorization of $f$) and $(Y,d_{Y})$ (a factorization of $g$) to produce a factorization, say $(Z,d_{Z})$, of $f+g$? \\
   Such details are not given in \cite{carqueville2016adjunctions}. These details are important to have a better understanding of the composition of $1-$morphisms in the bicategory of Landau Ginzburg models.  \\
  In this subsection, we explain how it is done. To that end, we need to prove lemma \ref{lemma (X tensor Y,D) is an obj of MF(f+g)} below, which actually shows how to produce a factorization of the sum of two polynomials from their respective factorizations, when the factorizations are defined as pairs $(X,d_{X})$ as earlier explained. \\
  Furthermore, using the variants of the Yoshino tensor product (one at a time), we also show how to produce a factorization of the sum of two polynomials from their respective factorizations, when the factorizations are defined as pairs $(X,d_{X})$.
  \\
Let $R$, $T$ and $S$ be commutative $K-$algebras with $f\in R=K[x]$ and $g\in S=K[y]$.
Let $(X,d_{X})$ be a factorization of $f \in K[x]$ and $(Y,d_{Y})$ be a factorization of $g \in K[y]$, where  we assume that $X$ is an $R-S-$bimodule and $Y$ is an $S-T-$bimodule. If we use the same notation for a map and the matrix representing it (precisely for $d_{X}$ and $d_{Y}$), then, we can use the Yoshino tensor product to define a new factorization $(X,d_{X}) \widehat{\otimes} (Y,d_{Y}):=(X\otimes Y, d_{X\otimes Y}=D)$ as follows:
$$(X\otimes Y)^{0}=(X^{0}\otimes Y^{0})\oplus (X^{1}\otimes Y^{1})\;\;\;\;\;(X\otimes Y)^{1}=(X^{0}\otimes Y^{1})\oplus (X^{1}\otimes Y^{0})$$
where the tensor product between the components of $X$ and $Y$ is taken over $S$. The matrix $D$ of the differential is obtained as follows:
Observe that looking at the matrices representing the Yoshino tensor product of matrix factorizations,
and mindful of the fact that $(\phi,\psi)$ and $(\phi',\psi')$ in that definition are actually represented here by $(d_{X}^{0},d_{X}^{1})$ and $(d_{Y}^{0},d_{Y}^{1})$, if we write $D=(D^{0},D^{1})$ the matrix representing the map $d_{X\otimes Y}$, then, we can write:
$$D^{1}= \begin{pmatrix}
d^{1}_{X}\otimes id_{Y^{1}} & -id_{X^{0}}\otimes d^{0}_{Y} \\
id_{X^{1}}\otimes d^{1}_{Y} & d^{0}_{X}\otimes id_{Y^{0}}

\end{pmatrix}\;\;\;\;\;\;D^{0}= \begin{pmatrix}
d^{0}_{X}\otimes id_{Y^{1}} & id_{X^{1}}\otimes d^{0}_{Y}\\
-id_{X^{0}}\otimes d^{1}_{Y} & d^{1}_{X}\otimes id_{Y^{0}}

\end{pmatrix}$$
where we used the same notation for a matrix and its map. In fact, for $i\in \{0,1\}$, $d_{Y}^{i}$ and $d_{X}^{i}$ are respectively the matrices of the maps $d_{Y}^{i}: Y^{i} \rightarrow Y^{i+1}$ and $d_{X}^{i}: X^{i} \rightarrow X^{i+1}$. Likewise $id_{Y^{i}}$ and $id_{X^{i}}$ are the matrices of the corresponding identity maps.\\
In the sequel, we will denote the units of $R$ and $S$ by $1_{R}$ and $1_{S}$ respectively.\\

\begin{lemma} \label{lemma (X tensor Y,D) is an obj of MF(f+g)}
$(X\otimes Y, D)$ as defined above, determines an object of
$MF(R\otimes S,f\otimes 1_{S} + 1_{R}\otimes g)$
\end{lemma}

\begin{proof}

We need to verify that $(X\otimes Y, D)$ is a factorization of $f\otimes 1_{S} + 1_{R}\otimes g \in R\otimes S.$
Let $h=f\otimes 1_{S} + 1_{R}\otimes g\in R\otimes S$.\\
Before proceeding with the proof of this lemma, recall what we saw immediately after definition \ref{finite rank matrix facto}, namely that
for a matrix factorization $X=X^{0}\oplus X^{1}$, its differential $d_{X}$ could be viewed as a matrix (that we still denote by $d_{X}$) and the condition to be fulfilled should then be $d^{2}_{X}=f\cdot Id$.
Now that we have the module $$X\otimes Y = (X^{0}\otimes Y^{0})\oplus (X^{1}\otimes Y^{1})\oplus(X^{0}\otimes Y^{1})\oplus (X^{1}\otimes Y^{0}).$$ In order to prove the lemma, we will need to prove that $D^{2}= h \cdot Id$.\\
 We know that: $$D^{2}=\begin{pmatrix}
 0 & D^{1} \\
D^{0} &  0
\end{pmatrix}\begin{pmatrix}
 0 & D^{1} \\
D^{0} &  0
\end{pmatrix}
= \begin{pmatrix}
  D^{1}D^{0}& 0 \\
   0  & D^{0}D^{1}
\end{pmatrix}.$$

In the sequel, we are going to implicitly use the mixed-product property (cf. Lemma 4.2.10 of \cite{horn2012matrix}) of the tensor product.
\begin{lemma}(cf. Lemma 4.2.10 of \cite{horn2012matrix})\\
 If $A,B,C$ and $D$ are matrices of such size that one can form the matrix products $AC$ and $BD$,
  then the product of two tensor products
yields another tensor product
$(A \otimes B)(C \otimes D) = AC \otimes BD$.
\end{lemma}
$$D^{1}D^{0}=  \begin{pmatrix}
d^{1}_{X}\otimes id_{Y^{1}} & -id_{X^{0}}\otimes d^{0}_{Y} \\
id_{X^{1}}\otimes d^{1}_{Y} & d^{0}_{X}\otimes id_{Y^{0}}

\end{pmatrix}\begin{pmatrix}
d^{0}_{X}\otimes id_{Y^{1}} & id_{X^{1}}\otimes d^{0}_{Y}\\
-id_{X^{0}}\otimes d^{1}_{Y} & d^{1}_{X}\otimes id_{Y^{0}}

\end{pmatrix}$$
%

$$= \begin{pmatrix}
d^{1}_{X}d^{0}_{X}\otimes id_{Y^{1}}id_{Y^{1}} + id_{X^{0}}id_{X^{0}}\otimes d^{0}_{Y}d^{1}_{Y} & d^{1}_{X}id_{X^{1}}\otimes id_{Y^{1}}d^{0}_{Y} - id_{X^{0}}d^{1}_{X}\otimes d^{0}_{Y}id_{Y^{0}} \\
id_{X^{1}}d^{0}_{X}\otimes d^{1}_{Y}id_{Y^{1}} - d^{0}_{X}id_{X^{0}}\otimes id_{Y^{0}}d^{1}_{Y} & id_{X^{1}}id_{X^{1}}\otimes d^{1}_{Y}d^{0}_{Y}
 + d^{0}_{X}d^{1}_{X}\otimes id_{Y^{0}}id_{Y^{0}}
\end{pmatrix}$$

$$= \begin{pmatrix}
f\cdot id_{X^{0}}\otimes id_{Y^{1}} + id_{X^{0}}\otimes g\cdot id_{Y^{1}}  & d^{1}_{X}\otimes d^{0}_{Y}- d^{1}_{X}\otimes d^{0}_{Y}\\
d^{0}_{X}\otimes d^{1}_{Y}- d^{0}_{X}\otimes d^{1}_{Y} & id_{X^{1}}\otimes g\cdot id_{Y^{0}} + f\cdot id_{X^{1}}\otimes id_{Y^{0}}

\end{pmatrix}$$

$$= \begin{pmatrix}
f\cdot id_{X^{0}}\otimes id_{Y^{1}} + id_{X^{0}}\otimes g\cdot id_{Y^{1}} & 0\\
0 & id_{X^{1}}\otimes g\cdot id_{Y^{0}} + f\cdot id_{X^{1}}\otimes id_{Y^{0}}
\end{pmatrix}$$
$$= \begin{pmatrix}
 (f\otimes 1_{S})\cdot (id_{X^{0}}\otimes id_{Y^{1}}) + (1_{R}\otimes g)\cdot(id_{X^{0}}\otimes id_{Y^{1}}) & 0\\
0 & (1_{R}\otimes g)\cdot(id_{X^{1}}\otimes id_{Y^{0}}) + (f\otimes 1_{S})\cdot (id_{X^{1}}\otimes id_{Y^{0}})
\end{pmatrix}$$
$$= \begin{pmatrix}
(f\otimes 1_{S} + 1_{R}\otimes g)\cdot(id_{X^{0}}\otimes id_{Y^{1}}) & 0\\
0 & (f\otimes 1_{S} + 1_{R}\otimes g)\cdot (id_{X^{1}}\otimes id_{Y^{0}})
\end{pmatrix}$$
$$= (f\otimes 1_{S} + 1_{R}\otimes g)\cdot\begin{pmatrix}
 id_{X^{0}\otimes Y^{1}} & 0\\
0 & id_{X^{1}\otimes Y^{0}}
\end{pmatrix}$$
$$= h\cdot\begin{pmatrix}
 id_{X^{0}\otimes Y^{1}} & 0\\
0 & id_{X^{1}\otimes Y^{0}}
\end{pmatrix}\,\,\,\,\,\,\,\,...\,\,\,\,\,\,\,\,\,(i)$$

$$D^{0}D^{1}=\begin{pmatrix}
d^{0}_{X}\otimes id_{Y^{1}} & id_{X^{1}}\otimes d^{0}_{Y}\\
-id_{X^{0}}\otimes d^{1}_{Y} & d^{1}_{X}\otimes id_{Y^{0}}

\end{pmatrix} \begin{pmatrix}
d^{1}_{X}\otimes id_{Y^{1}} & -id_{X^{0}}\otimes d^{0}_{Y} \\
id_{X^{1}}\otimes d^{1}_{Y} & d^{0}_{X}\otimes id_{Y^{0}}

\end{pmatrix}$$
%

$$= \begin{pmatrix}
d^{0}_{X}d^{1}_{X}\otimes id_{Y^{1}}id_{Y^{1}} + id_{X^{1}}id_{X^{1}}\otimes d^{0}_{Y}d^{1}_{Y} & -d^{0}_{X}id_{X^{0}}\otimes id_{Y^{1}}d^{0}_{Y} + id_{X^{1}}d^{0}_{X}\otimes d^{0}_{Y}id_{Y^{0}} \\
-id_{X^{0}}d^{1}_{X}\otimes d^{1}_{Y}id_{Y^{1}} + d^{1}_{X}id_{X^{1}}\otimes id_{Y^{0}}d^{1}_{Y} & id_{X^{0}}id_{X^{0}}\otimes d^{1}_{Y}d^{0}_{Y}
 + d^{1}_{X}d^{0}_{X}\otimes id_{Y^{0}}id_{Y^{0}}
\end{pmatrix}$$

$$= \begin{pmatrix}
f\cdot id_{X^{1}}\otimes id_{Y^{1}} + id_{X^{1}}\otimes g\cdot id_{Y^{1}}  & -d^{0}_{X}\otimes d^{0}_{Y}+ d^{0}_{X}\otimes d^{0}_{Y}\\
-d^{1}_{X}\otimes d^{1}_{Y}+ d^{1}_{X}\otimes d^{1}_{Y} & id_{X^{0}}\otimes g\cdot id_{Y^{0}} + f\cdot id_{X^{0}}\otimes id_{Y^{0}}

\end{pmatrix}$$

$$= \begin{pmatrix}
f\cdot id_{X^{1}}\otimes id_{Y^{1}} + id_{X^{1}}\otimes g\cdot id_{Y^{1}} & 0\\
0 & id_{X^{0}}\otimes g\cdot id_{Y^{0}} + f\cdot id_{X^{0}}\otimes id_{Y^{0}}
\end{pmatrix}$$
$$= \begin{pmatrix}
 (f\otimes 1_{S})\cdot (id_{X^{1}}\otimes id_{Y^{1}}) + (1_{R}\otimes g)\cdot(id_{X^{1}}\otimes id_{Y^{1}}) & 0\\
0 & (1_{R}\otimes g)\cdot(id_{X^{0}}\otimes id_{Y^{0}}) + (f\otimes 1_{S})\cdot (id_{X^{0}}\otimes id_{Y^{0}})
\end{pmatrix}$$
$$= \begin{pmatrix}
(f\otimes 1_{S} + 1_{R}\otimes g)\cdot(id_{X^{1}}\otimes id_{Y^{1}}) & 0\\
0 & (f\otimes 1_{S} + 1_{R}\otimes g)\cdot (id_{X^{0}}\otimes id_{Y^{0}})
\end{pmatrix}$$
$$= (f\otimes 1_{S} + 1_{R}\otimes g)\cdot\begin{pmatrix}
 id_{X^{1}\otimes Y^{1}} & 0\\
0 & id_{X^{0}\otimes Y^{0}}
\end{pmatrix}$$
$$= h\cdot\begin{pmatrix}
 id_{X^{1}\otimes Y^{1}} & 0\\
0 & id_{X^{0}\otimes Y^{0}}
\end{pmatrix}\,\,\,\,\,\,\,\,...\,\,\,\,\,\,\,\,\,(ii)$$

Hence, putting $(i)$ and $(ii)$ together, we get the desired result, that is;

$$D^{2}=\begin{pmatrix}
 0 & D^{1} \\
D^{0} &  0
\end{pmatrix}\begin{pmatrix}
 0 & D^{1} \\
D^{0} &  0
\end{pmatrix}
= \begin{pmatrix}
  D^{1}D^{0}& 0 \\
   0  & D^{0}D^{1}
\end{pmatrix}
= h\cdot\begin{pmatrix}
 id_{X^{0}\otimes Y^{1}} & 0 & 0& 0\\
0 & id_{X^{1}\otimes Y^{0}} &0 &0 \\
0 &0 & id_{X^{1}\otimes Y^{1}}\\
0 & 0& 0 & id_{X^{0}\otimes Y^{0}}
\end{pmatrix}=h\cdot Id$$

\end{proof}

Keeping the notations above, we can use any of the three variants of the Yoshino tensor product to define new factorizations which differ at the level of their differentials. Thus for the first variant we have: $(X,d_{X})\otimes' (Y,d_{Y}):=(X\otimes ' Y, d_{X\otimes' Y}=\underline{D})$ and
$$(X\otimes ' Y)^{0}=(X^{0}\otimes Y^{0})\oplus (X^{1}\otimes Y^{1})\;\;\;\;\;(X\otimes ' Y)^{1}=(X^{0}\otimes Y^{1})\oplus (X^{1}\otimes Y^{0})$$
where the tensor product between the components of $X$ and $Y$ is taken over $S$. The matrix $\underline{D}=(\underline{D}^{0},\underline{D}^{1})$ of the differential is obtained as follows:\\
$\underline{D}^{0}$ is obtained by "rotating" $D^{0}$ anticlockwise once, while $\underline{D}^{1}$ is obtained by "rotating" $D^{1}$ clockwise once. Thus:
$$\underline{D}^{1}= \begin{pmatrix}
id_{X^{1}}\otimes d^{1}_{Y} & d^{1}_{X}\otimes id_{Y^{1}} \\
d^{0}_{X}\otimes id_{Y^{0}} & -id_{X^{0}}\otimes d^{0}_{Y}

\end{pmatrix}\;\;\;\;\;\;\underline{D}^{0}= \begin{pmatrix}
id_{X^{1}}\otimes d^{0}_{Y} & d^{1}_{X}\otimes id_{Y^{0}}\\
d^{0}_{X}\otimes id_{Y^{1}} & -id_{X^{0}}\otimes d^{1}_{Y}

\end{pmatrix}$$
where we used the same notation for a matrix and its map as done in the case of the Yoshino tensor product above. \\
Analogously, for the second variant, the matrix $\underline{\underline{D}}=(\underline{\underline{D}}^{0},\underline{\underline{D}}^{1})$ of the differential is obtained as follows:\\
$\underline{\underline{D}}^{0}$ is obtained by "rotating" $D^{0}$ anticlockwise twice, while $\underline{\underline{D}}^{1}$ is obtained by "rotating" $D^{1}$ clockwise twice. Thus:
$$\underline{\underline{D}}^{1}= \begin{pmatrix}
d^{0}_{X}\otimes id_{Y^{0}} & id_{X^{1}}\otimes d^{1}_{Y} \\
-id_{X^{0}}\otimes d^{0}_{Y} & d^{1}_{X}\otimes id_{Y^{1}}

\end{pmatrix}\;\;\;\;\;\;\underline{\underline{D}}^{0}= \begin{pmatrix}
d^{1}_{X}\otimes id_{Y^{0}} & -id_{X^{0}}\otimes d^{1}_{Y}\\
id_{X^{1}}\otimes d^{0}_{Y} & d^{0}_{X}\otimes id_{Y^{1}}

\end{pmatrix}$$
where we used the same notation for a matrix and its map.\\
We omit the case of the third variant, because it is analogously obtained. In fact, to find the matrix of the differential, it suffices to "rotate" three times $D^{0}$ anticlockwise and $D^{1}$ clockwise.

An analogue of Lemma \ref{lemma (X tensor Y,D) is an obj of MF(f+g)} can now be stated for each of the three variants. We do it only for the second variant since the other ones are proved similarly.
\begin{lemma} \label{lemma (X tensor'' Y,D) is an obj of MF(f+g)}
$(X\otimes '' Y, D)$ as defined above, determines an object of
$MF(R\otimes S,f\otimes 1_{S} + 1_{R}\otimes g)$
\end{lemma}

\begin{proof} (A sketch)

We need to prove that $\underline{\underline{D}}^{2}= h \cdot Id$, where $h=f\otimes 1_{S} + 1_{R}\otimes g\in R\otimes S$.\\
 We know that: $$\underline{\underline{D}}^{2}=\begin{pmatrix}
 0 & \underline{\underline{D}}^{1} \\
\underline{\underline{D}}^{0} &  0
\end{pmatrix}\begin{pmatrix}
 0 & \underline{\underline{D}}^{1} \\
\underline{\underline{D}}^{0} &  0
\end{pmatrix}
= \begin{pmatrix}
  \underline{\underline{D}}^{1}\underline{\underline{D}}^{0}& 0 \\
   0  & \underline{\underline{D}}^{0}\underline{\underline{D}}^{1}
\end{pmatrix}.$$

In the sequel, we are going to implicitly use the mixed-product property (cf. Lemma 4.2.10 in \cite{horn2012matrix}) of the tensor product.
$$\underline{\underline{D}}^{0}\underline{\underline{D}}^{1}=\begin{pmatrix}
d^{1}_{X}\otimes id_{Y^{0}} & -id_{X^{0}}\otimes d^{1}_{Y}\\
id_{X^{1}}\otimes d^{0}_{Y} & d^{0}_{X}\otimes id_{Y^{1}}

\end{pmatrix}\begin{pmatrix}
d^{0}_{X}\otimes id_{Y^{0}} & id_{X^{1}}\otimes d^{1}_{Y} \\
-id_{X^{0}}\otimes d^{0}_{Y} & d^{1}_{X}\otimes id_{Y^{1}}

\end{pmatrix}$$
$$= \begin{pmatrix}
d^{1}_{X}d^{0}_{X}\otimes id_{Y^{0}} + id_{X^{0}}\otimes d^{1}_{Y}d^{0}_{Y} & 0 \\
0 & id_{X^{1}}\otimes d^{0}_{Y}d^{1}_{Y}
 + d^{0}_{X}d^{1}_{X}\otimes id_{Y^{1}}
\end{pmatrix}$$

$$= \begin{pmatrix}
f\cdot id_{X^{0}}\otimes id_{Y^{0}} + id_{X^{0}}\otimes g\cdot id_{Y^{0}}  & 0\\
0 & id_{X^{1}}\otimes g\cdot id_{Y^{1}} + f\cdot id_{X^{1}}\otimes id_{Y^{1}}

\end{pmatrix}$$

$$= (f\otimes 1_{S} + 1_{R}\otimes g)\cdot\begin{pmatrix}
 id_{X^{0}\otimes Y^{0}} & 0\\
0 & id_{X^{1}\otimes Y^{1}}
\end{pmatrix}$$
$$= h\cdot\begin{pmatrix}
 id_{X^{0}\otimes Y^{0}} & 0\\
0 & id_{X^{1}\otimes Y^{1}}
\end{pmatrix}\,\,\,\,\,\,\,\,...\,\,\,\,\,\,\,\,\,(i)$$

$$\underline{\underline{D}}^{1}\underline{\underline{D}}^{0}=\begin{pmatrix}
d^{0}_{X}d^{1}_{X}\otimes id_{Y^{0}} + id_{X^{1}}\otimes d^{1}_{Y}d^{0}_{Y} & 0\\
0 & id_{X^{0}}\otimes d^{0}_{Y}d^{1}_{Y} + d^{1}_{X}d^{0}_{X}\otimes id_{Y^{1}}

\end{pmatrix} $$

$$= \begin{pmatrix}
f\cdot id_{X^{1}}\otimes id_{Y^{0}} + id_{X^{1}}\otimes g\cdot id_{Y^{0}} & 0\\
0 & id_{X^{0}}\otimes g\cdot id_{Y^{1}} + f\cdot id_{X^{0}}\otimes id_{Y^{1}}
\end{pmatrix}$$

$$= h\cdot\begin{pmatrix}
 id_{X^{1}\otimes Y^{0}} & 0\\
0 & id_{X^{0}\otimes Y^{1}}
\end{pmatrix}\,\,\,\,\,\,\,\,...\,\,\,\,\,\,\,\,\,(ii)$$

Hence, putting $(i)$ and $(ii)$ together, we get the desired result, that is;

$$\underline{\underline{D}}^{2}=\begin{pmatrix}
 0 & \underline{\underline{D}}^{1} \\
\underline{\underline{D}}^{0} &  0
\end{pmatrix}\begin{pmatrix}
 0 & \underline{\underline{D}}^{1} \\
\underline{\underline{D}}^{0} &  0
\end{pmatrix}
= \begin{pmatrix}
  \underline{\underline{D}}^{1}\underline{\underline{D}}^{0}& 0 \\
   0  & \underline{\underline{D}}^{0}\underline{\underline{D}}^{1}
\end{pmatrix}
= h\cdot\begin{pmatrix}
 id_{X^{1}\otimes Y^{0}} & 0 & 0& 0\\
0 & id_{X^{0}\otimes Y^{1}} &0 &0 \\
0 &0 & id_{X^{0}\otimes Y^{0}}\\
0 & 0& 0 & id_{X^{1}\otimes Y^{1}}
\end{pmatrix}=h\cdot Id$$

\end{proof}
The case of the first and third variant are proved similarly. But it is important to note that in these two cases one needs to observe (thanks to the discussion preceding lemma \ref{lemma (X tensor Y,D) is an obj of MF(f+g)}) that the identity matrices $id_{X^{0}}$ and $id_{X^{1}}$ are of the same size. The same holds for $id_{Y^{0}}$ and $id_{Y^{1}}$.\\
We will soon allude to this subsection \ref{subsec:tensor prod} when defining the composition of $1-$morphisms in $\mathcal{LG}_{K}$ which are matrix factorizations. In fact, we will observe that instead of using the Yoshino tensor product as done in \cite{carqueville2016adjunctions} to define the composition of two matrix factorizations, we could also use any of its three variants as discussed above. This will naturally give rise to three variants of $\mathcal{LG}_{K}$.\\

The following properties are proved in \cite{yoshino1998tensor} and the tensor product here is the Yoshino tensor product.
\begin{lemma} (section 2 of \cite{yoshino1998tensor})
  \begin{itemize}
    \item If $h:(X,d)\longrightarrow (Y,d')$ is a morphism of factorizations of $f\in R$ and $(Z,d")$ is a factorization of $g\in S$, then there is an evident morphism of factorizations $h\otimes Z$. Likewise, if
        $r:(T,d_{1})\longrightarrow (L,d_{1}')$ is a morphism of factorizations of $g\in S$ and $(Z_{1},d_{1}")$ is a factorization of $f\in R$, then there is an evident morphism of factorizations $Z_{1}\otimes r$.
    \item Let $(X,d)$ be a factorization of $f$ and $(Y,d')$ be a factorization of $g$. Then there is an isomorphism $X\otimes Y \cong Y\otimes X$.
    \item There is an isomorphism of factorizations $(X\otimes Y)\otimes Z\cong X\otimes (Y\otimes Z)$
    \item The category F($R,f$) has biproducts which are obtained in the evident way, and tensor distributes over this biproduct:
        $$X\otimes (Y\oplus Z)\cong (X\otimes Y)\oplus (X\otimes Z)$$
  \end{itemize}
\end{lemma}

It is time to review the bicategory $\mathcal{LG}_{K}$ constructed in \cite{carqueville2016adjunctions}.
As mentioned at the beginning of this paper, in our presentation of $\mathcal{LG}_{K}$, objects are polynomials without the restrictions imposed on them in section 2.2 of \cite{carqueville2016adjunctions}. This poses no problem to what we want to do.

\section{A review of the construction of the bicategory $\mathcal{LG}_{K}$} \label{bicategory LG}
We proceed in steps. First, we recall the definition of a bicategory. Next, we construct a structure that has all the ingredients of a bicategory except for the existence of identity one-cells. Thereafter, we develop a sophisticated machinery to elucidate the intricate construction of the entity that has to act as unit in the bicategory $\mathcal{LG}_{K}$. Finally, we prove that there is no direct inverse for the unitors (i.e., the right and left identities, see \cite{carqueville2016adjunctions}) thereby explaining why in \cite{carqueville2016adjunctions}, their construction is done at the level of homotopy.

\begin{definition} \label{bicat} \cite{barbosa2003brief}\\
A \textbf{bicategory} $\mathcal{B}$ is made up of the following data:
\begin{enumerate}
\item A class of objects $A, B, C,... $
\item For each pair $\langle A, B \rangle$ of objects, a small (hom-)category $\mathcal{B}(A, B)$ with
arrows $p,q,r,...$ from $A$ to $B$ as objects and arrows $\alpha,\beta,\gamma,...$ between them, referred to as 2-cells.
The arrows $p,q,r,...$ are called 1-cells. \\An example of a 2-cell is denoted as follows: $\alpha: p\Longrightarrow q$.
Composition in $\mathcal{B}(A, B)$ is denoted by a dot "$\cdot $" and the identity on $p$ for each $p: A\longrightarrow B$, by $ 1_{p}: p \Longrightarrow p $.

\item For each object $A$, a functor that returns the identity on $A$.

$ I_{A}: 1 \rightarrow \mathcal{B}(A, A)$, where $1$ stands for the final object in the category Cat of small categories.

If we write $1=\{\ast\}$, that is the unique object of $1$ is denoted $\ast$, then $\mathbf{I}_{A} (\ast)$ is
the 1-cell $1_{A}: A\longrightarrow A$.

\item For each triple $\langle A, B, C \rangle$ of objects, a composition law given by a functor \\
$\star_{A,B,C}: \mathcal{B}(A, B)\times \mathcal{B}(B, C)\longrightarrow \mathcal{B}(A, C)$.

\item For each quadruple $\langle A,B,C,D \rangle$ of objects of $\mathcal{B}$, a natural isomorphism
%

 $$a_{A,B,C,D}: \star_{A,B,D}\circ (Id \times \star_{B,C,D})\Longrightarrow \star_{A,C,D}\circ (\star_{A,B,C}\times Id)$$ where Id is the identity functor.\\
The components of this isomorphism are the invertible 2-cells defined as follows: For $p:A\longrightarrow B$, $q:B\longrightarrow C$, $r:C\longrightarrow D$ and $t:D\longrightarrow E$,  we have
$$ \xymatrix { a_{p,q,r}:= a_{A,B,C,D}(p,q,r): \star_{A,B,D}\circ (Id \times \star_{B,C,D})(p,q,r)
\ar @{=>}[r]^(.63){\cong} &\star_{A,C,D} \circ (\star_{A,B,C}\times Id)(p,q,r) }$$
$$i.e.,\; \xymatrix { a_{p,q,r}: p\star(q\star r) \ar @{=>}[r]^(.57){\cong} & (p\star q)\star r }$$
\item For each pair $\langle A,B\rangle$ of objects of $\mathcal{B}$, two natural isomorphisms $l(A,B)$ and
$r(A,B)$, called \textit{left} and \textit{right} isomorphisms or identities or unit laws (or unitors in \cite{carqueville2016adjunctions}). They are given respectively by:\\
$$l_{A,B}: \star_{A,B,B}\circ (Id\times I_{B})\Longrightarrow Id$$
$$ r_{A,B}: \star_{A,A,B}\circ (I_{A}\times Id)\Longrightarrow Id$$
%
%
%
The components of these isomorphisms are the invertible 2-cells defined as follows: For $p:A\longrightarrow B$, we have
%
%
$$\; \xymatrix {r_{p}: 1_{A}\star p\ar @{=>}[r]^(.625){\cong} & p} $$
%
%
%
and
$$\xymatrix {l_{p}:  p \star 1_{B}\ar @{=>}[r]^(.63){\cong} & p} $$
The families of natural isomorphisms $a_{A,B,C,D}$ , $r_{A,B}$ and $l_{A,B}$ are furthermore required to satisfy the associativity and identity coherence laws (cf. definition 1.1 \cite{fomatati2019multiplicative}).

\end{enumerate}
\end{definition}

The following lemma is useful.
\begin{lemma}\label{free mod isos} \cite{conrad2016tensor}
\begin{enumerate}
  \item Let $K[x_{1},x_{2},\cdots,x_{n}]$ and $K[x'_{1},x'_{2},\cdots,x'_{m}]$ be free $K-$modules with respective bases $\{e_{i}\}^{n}_{i=1}$ and $\{e'_{j}\}^{m}_{j=1}$. Then $\{e_{i}\otimes e'_{j}\}_{i=1,\cdots n; j=1,\cdots m}$ is a basis of $K[x_{1},x_{2},\cdots,x_{n}]\otimes K[x'_{1},x'_{2},\cdots,x'_{m}]$.
  \item The $K-$modules $K[x_{1},x_{2},\cdots,x_{n}]\otimes K[x'_{1},x'_{2},\cdots,x'_{m}]$ and $K[x_{1},x_{2},\cdots,x_{n},x'_{1},x'_{2},\cdots,x'_{m}]$ are isomorphic as $K-$modules.

  \item If we let $x$ stand for $x_{1},x_{2},\cdots,x_{n}$ and $x^{(l)}$ stand for $x^{(l)}_{1},x^{(l)}_{2},\cdots, x^{(l)}_{n_{l}}$, where $n,l,n_{l}\in \mathbb{N}$, $n=n_{0}$, and $\{x^{(l)}\}$ means $x$ with $l$ primes, e.g $x^{(2)}=x'',\;x^{(0)}:=x$.\\
       Then more generally, we have:
      $$K[x,x^{(1)},\cdots, x^{(l)}]\cong \bigotimes_{p=0,\cdots,l} K[x^{(p)}].$$
\end{enumerate}

\end{lemma}
%
%
%

\subsection{Bicategorical structure}
Here, we introduce the notion of $B-$category (cf. definition \ref{defn B-category}) in order to facilitate the discussion on the bicategorical structure of $\mathcal{LG}_{K}$.
\begin{definition}$B-$category  \label{defn B-category} \\
  A \textbf{$B-$category} $\mathcal{B}$ is an algebraic structure that has all the structure of a bicategory as defined in definition \ref{bicat} except for the existence of identity one-cells.
\end{definition}
That is, $\mathcal{B}$ is made up of the data given in definition \ref{bicat} except point $3$, and we should also stress in point $2$ that identity one-cells are not required to exist.\\
Before we continue, we make the following useful remark about the homotopy category of matrix factorizations ($HMF(R,f)$) and one of its interesting full subcategories; namely the homotopy category of finite rank matrix factorizations that is denoted by $hmf(R,f)$.
\begin{remark} (p.9 of \cite{carqueville2016adjunctions})
$HMF(R,f)$ is idempotent complete (\cite{bokstedt1993homotopy}, \cite{neeman2001triangulated}).
Since we work with polynomials rather than power series, $hmf(R,f)$ is not necessarily idempotent complete \cite{keller2011two}. The idempotent closure of $hmf(R,f)$ denoted by $hmf(R,f)^{\omega}$ is a full subcategory of $HMF(R,f)$ whose objects are those matrix factorizations which are direct summands of finite-rank matrix factorizations in the homotopy category. Moreover, $hmf(R,f)^{\omega}$ is an idempotent complete category.\\
As explained in \cite{carqueville2016adjunctions} (p.9), taking the idempotent completion is necessary because the composition of $1$-morphisms in $\mathcal{LG}_{K}$ results in matrix factorisations which, while not finite-rank, are summands in the homotopy category of something finite-rank. There are two natural ways to resolve this: work throughout with power series rings and completed tensor products, or work
with idempotent completions.

\end{remark}

We construct a $B-$category which we call $B-Fac$. The objects of $B-Fac$ are polynomials $f$ denoted by pairs $(R,f)$ where $f\in R=K[x]$. Let ($R=K[x], f$) and ($S=K[z], g$) be elements of $B-Fac$. We then define the small category $B-Fac((R,f),(S,g))$ as follows:
$$B-Fac((R,f),(S,g)):=hmf(R\otimes S, 1_{R} \otimes g - f\otimes 1_{S} )^{\omega}=hmf(K[x,z] , g - f )^{\omega}$$
viz. a 1-morphism between two polynomials $f$ and $g$ is a matrix factorization of $g-f$. (Recall that $hmf(K[x,z] , g - f )^{\omega}$ is a subcategory of the category of matrix factorizations modulo homotopy denoted by $HMF(K[x,z], g - f )$).\\
Then given two composable $1-$cells $X\in B-Fac((R,f),(S,g))$ and \\ $Y\in B-Fac((S,g),(T,h))$, we can define their composition using Yoshino's tensor product as it is done in \cite{carqueville2016adjunctions}. In fact, by remark 2.1.8 on p.29 of \cite{camacho2015matrix}, $Y\otimes_{S} X$ is a free module of infinite rank over $R\otimes_{K} T$.
 However, the argument of Section 12 of \cite{dyckerhoff2013pushing} shows that it is naturally isomorphic to a direct summand in the homotopy category of something finite-rank. Thus, we may define $Y\circ X:=Y\otimes_{S} X$ $\in\,hmf(R\otimes_{K}T, 1_{R}\otimes h - f\otimes 1_{T})^{\omega}= B-Fac((R,f),(T,h))$.\\
We can also define their composition using variants of Yoshino's tensor product as discussed in subsection \ref{subsec:tensor prod}. \\
So,
$(Y, d_{Y})\circ (X, d_{X})$ could be chosen to be any of the following four tensor products: $(Y, d_{Y})\widehat{\otimes} (X, d_{X})= (Y\otimes_{S} X, d_{Y\widehat{\otimes} X})$, $(Y, d_{Y})\widehat{\otimes}' (X, d_{X})= (Y\otimes_{S} X, d_{Y\widehat{\otimes}' X})$, $(Y, d_{Y})\widehat{\otimes}'' (X, d_{X})= (Y\otimes_{S} X, d_{Y\widehat{\otimes}'' X})$ and $(Y, d_{Y})\widehat{\otimes}''' (X, d_{X})= (Y\otimes_{S} X, d_{Y\widehat{\otimes}''' X})$
 which are all objects of  $HMF(R\otimes_{K}T, 1_{R}\otimes h - f\otimes 1_{T})$ and are $\mathbb{Z}_{2}-$graded modules, where
$$(Y\otimes_{S} X)^{0}=(Y^{0}\otimes_{S} X^{0})\oplus (Y^{1}\otimes_{S} X^{1})\;\;and \;\; (Y\otimes_{S} X)^{1}=(Y^{0}\otimes_{S} X^{1})\oplus (Y^{1}\otimes_{S} X^{0})$$
The differentials for these tensor products (Yoshino's product and its variants) are different and they are defined in subsection \ref{subsec:tensor prod} above.\\
\textbf{N.B.} The above different ways to define composition of $1-$morphisms give rise to three variants of the bicategory $\mathcal{LG}_{K}$. In fact, each variant of the Yoshino tensor product gives rise to a variant of
$\mathcal{LG}_{K}$. In the presentation of the construction of this bicategory below, one can observe that if any of the foregoing variants is used in the construction of $\mathcal{LG}_{K}$, then we will obtain virtually the same bicategory. The result is simply a variant of $\mathcal{LG}_{K}$. So, we have four versions of $\mathcal{LG}_{K}$ which only differ at the level of the definition of the composition of $1-$cells.
\\

It is now time to define the tensor product of morphisms of matrix factorizations. Let $X_{1}$, $X_{2}$ be
objects of $B-Fac((R,f),(S,g))$ and $Y_{1}$, $Y_{2}$ be
objects of $B-Fac((S,g), (T,h))$. Let $\alpha : X_{1} \rightarrow X_{2}$ and $\beta : Y_{1} \rightarrow Y_{2}$ be two morphisms, then we define their tensor product in the obvious way $\beta \otimes \alpha : Y_{1} \otimes X_{1} \rightarrow Y_{2}\otimes X_{2}$ in $B-Fac((R,f),(T,h))$.\\
With the above data, the composition (bi-)functor is entirely determined in our B-category:

$\star_{(R,f),(S,g),(T,h)}:B-Fac((R,f),(S,g))$ x  $B-Fac((S,g), (T,h))\rightarrow B-Fac((R,f),(T,h))$, $(X, Y)\mapsto Y\otimes_{S} X$.\\

The definition of the associativity morphism is easy to state. In fact, for $X\in B-Fac((R,f),(S,g))$, $Y\in B-Fac((S,g),(T,h))$ and $Z\in B-Fac((T,h), (P,r))$, the associator is the 2-isomorphism $$a_{Z,Y,X}: Z\otimes (Y\otimes X)\rightarrow (Z\otimes Y)\otimes X$$ given by the usual formula $$z\otimes (y\otimes x)\rightarrow (z\otimes y)\otimes x$$
where $x\in X$, $y\in Y$ and $z\in Z$.


\begin{lemma}
The structure $B-Fac$ is a B-category.
\end{lemma}

  This lemma is proved by showing that all the conditions required to obtain a bicategory are met except for the existence of units and the left and right unit actions. The proof though not difficult is lengthy and we would rather focus more on the units construction which is one of our main points of interest in the review of the  bicategory of Landau-Ginzburg models.\\

The rest of the paper is reserved to the construction of the unit $1-$morphisms of $\mathcal{LG}_{K}$ and to the proof that the right and left unit maps are natural. We will prove that these unit actions do not possess direct inverses thereby accounting for the fact that their inverses are found only up to homotopy in \cite{carqueville2016adjunctions}.\\

In order to obtain the bicategory $\mathcal{LG}_{K}$ from the B-category $B-Fac$, we need to define the unit $1-$morphisms, the left and the right unitors.
\subsection{Unit $1-$morphisms in $\mathcal{LG}_{K}$} \label{subsec: unit 1-morphisms in LGk}
In this subsection, we will construct the identity 1-cells including all the intricacies involved. Recall that $R=K[x_{1},x_{2},\cdots,x_{n}]$. \\

From lemma \ref{free mod isos}, we have in particular that $$R\otimes_{K} R \cong K[x_{1},x_{2},\cdots,x_{n},x'_{1},x'_{2},\cdots,x'_{n}].$$
where $x_{i}=x_{i}\otimes 1$ and $x^{'}_{i}= 1 \otimes x_{i}$.
\\
The subscript $"K"$ in $\otimes_{K}$ will be very often omitted for ease of notation.
We need an object $\Delta_{f}: (R,f) \rightarrow (R,f)$ in
$hmf(R\otimes R,  f\otimes id - id\otimes f)$ or equivalently,\\
$hmf(K[x_{1},x_{2},\cdots,x_{n},x'_{1},x'_{2},\cdots,x'_{n}], h(x,x')),$ where $h(x,x')= f(x)- f(x'),$ where $id$ stands for $1_{R}$. In the remainder of this paper, $id$ denotes the identity in the ring under consideration.
\\
Recall (cf. section 5.5 of \cite{smith2011introduction}): The exterior algebra $\bigwedge(V)$ of a vector space $V$ over a field K is defined as the quotient algebra of the tensor algebra; $T(V)=\bigoplus_{i=1}^{\infty}T^{i}(V)=K\bigoplus V \bigoplus (V\otimes V)\bigoplus (V\otimes V\otimes V)\bigoplus \cdots ,$ by the two-sided ideal $I$ generated by all elements of the form $x \otimes x$ for $x \in V$. Symbolically, $\bigwedge(V)=T(V)/I$. The exterior product $\wedge$ of two elements of $\bigwedge(V)$ is the product induced by the tensor product $\otimes$ of $T(V)$. That is, if

     $\pi: T ( V )\rightarrow \bigwedge( V ) = T ( V ) / I $

is the canonical surjection, and if a and b are in $\bigwedge(V)$, then there are $\alpha$  and  $\beta$ in $T(V)$ such that $a = \pi ( \alpha )$ and $b = \pi ( \beta )$ and $a\wedge b=\pi(\alpha \otimes \beta)$.
 Let $\theta_{1},\theta_{2},\cdots, \theta_{n}$ be formal symbols that is, they are linearly independent by definition. We consider the $R\otimes R-$module:
 $$\Delta_{f}=\bigwedge\{\bigoplus_{i=1}^{n}(R\otimes R)\theta_{i}\}$$

 This is an exterior algebra generated by $n$ anti-commuting variables, the $\theta_{i}s$
 modulo the relations that the $\theta$'s anti-commute, that is $\theta_{i}\wedge \theta_{j}=-\theta_{j}\wedge \theta_{i}$. Typically, we will omit the wedge product and write for instance $\theta_{i}\wedge \theta_{j}$ simply as $\theta_{i}\theta_{j}$. Here, the "$\wedge$" product is taken over K just like the tensor product.
 A typical element of $\Delta_{f}$ is $(r\otimes r')\theta_{i_{1}}\theta_{i_{2}}\cdots \theta_{i_{k}}$ or equivalently $h(x_{1},x_{2},\cdots,x_{n},x'_{1},x'_{2},\cdots,x'_{n})\theta_{i_{1}}\theta_{i_{2}}\cdots \theta_{i_{k}}$ where $i_{1},\cdots, i_{k}\in \{1,\cdots,n\}$ and $h(x_{1},x_{2},\cdots,x_{n},x'_{1},x'_{2},\cdots,x'_{n})$ $\in$
 $K[x_{1},x_{2},\cdots,x_{n},x'_{1},x'_{2},\cdots,x'_{n}]$.\\
 $\Delta_{f}$ as an algebra is finitely generated by the set of formal symbols $\{\theta_{1},\cdots, \theta_{n}\}$. \\
$\Delta_{f}$ as an $R\otimes R-$module is generated by the set containing the empty list and all products of the form $\theta_{i_{1}}\dots \theta_{i_{k}}$ where $i_{1},\cdots, i_{k} \in \{1,\cdots,n\}$. This set of generators must be finite because by definition of $\Delta_{f}$, $\theta_{i_{p}} \theta_{i_{q}}=-\theta_{i_{p}} \theta_{i_{q}}$ and $\theta_{i_{p}} \theta_{i_{q}}=0$ if $p=q$ (${i_{p}}, {i_{q}} \in \{1,\cdots,n\}$). Moreover we know that since $\{\theta_{1},\dots, \theta_{n}\}$ is a set of formal symbols, it follows that the generating set of the $R\otimes R-$module $\Delta_{f}$ will be linearly independent and hence will form a finite basis. Whence, $\Delta_{f}$ is a finite rank $R\otimes R-$module.
 In fact:\\
  The action of $R\otimes R$ is the obvious one. That is, for example, if $r_{1}\otimes r_{2}\in R\otimes R$, its left action on an element $(r\otimes r')\theta_{i_{1}}\theta_{i_{2}}\cdots \theta_{i_{k}}$ of $\Delta_{f}$
  simply yields $(r_{1}\otimes r_{2})(r\otimes r')\theta_{i_{1}}\theta_{i_{2}}\cdots \theta_{i_{k}}$.
  The right action would yield $(r\otimes r')(r_{1}\otimes r_{2})\theta_{i_{1}}\theta_{i_{2}}\cdots \theta_{i_{k}}$ which is in fact the same as the left action since $R$ is commutative.

  $\Delta_{f}$ is endowed with the $\mathbb{Z}_{2}-$grading given by $\theta-$degree (where deg$\theta_{i}=1$ for each $i$). Thus $deg\theta_{i}^{2}=0\;and\;deg\theta_{i}\theta_{j}=0$.
 \\
  Next, we define the differential as follows:\\
  $$d: \Delta_{f}\longrightarrow \Delta_{f}$$
  $$d(-)=\sum_{i=1}^{n}[(x_{i}-x_{i}^{'})\theta_{i}^{\ast}(-)+\partial_{i}(f)\theta_{i}\wedge (-)]\,\,\,\,\,\,\,\,\,\,\,\,\,\,\,\,\cdots\,\,\,\natural$$
Where $\theta_{i}^{\ast}$ is the unique derivation extending the map $\theta_{i}^{\ast}(\theta_{j})=\delta_{ij}$ and as mentioned, it acts on an element $\theta_{i_{1}}\theta_{i_{2}}\cdots \theta_{i_{k}}$ of the
exterior algebra by the Leibniz rule with Koszul signs. Here, we elucidate what this means. In fact,
$$
\theta_{i}^{\ast}(\theta_{j_{1}}\theta_{j_{2}}\cdots \theta_{j_{k}})=
\begin{cases}
0\,for\,i\neq j,\,\forall j\in \{{j_{1}},{j_{2}},\cdots,{j_{k}}\}\\
(-1)^{p+1}\theta_{j_{1}}\theta_{j_{2}}\cdots \hat{\theta_{i}}\cdots \theta_{j_{k}}\,\,otherwise
\end{cases}
$$

where
 $\hat{\theta_{i}}$ signifies
that $\theta_{i}$ has been removed, and $p$ is the position of $\theta_{i}$ in $\theta_{j_{1}}\theta_{j_{2}}\cdots \theta_{i}\cdots \theta_{j_{k}}$\\

N.B. As earlier mentioned (\cite{carqueville2016adjunctions}), $\theta_{i}^{\ast}$ acts on an element $\theta_{i_{1}}\cdots\theta_{i_{k}}$ of the exterior algebra by the Leibniz rule with Koszul sign.
So, $\theta_{i}^{\ast}(\theta_{j}\theta_{k})=\theta_{i}^{\ast}(\theta_{j})\theta_{k}+(-1)^{\mid\theta_{j}\mid}\theta_{j}\theta_{i}^{\ast}(\theta_{k})$. Clearly, $\theta_{i}^{\ast}(\theta_{j}\theta_{k})=\theta_{i}^{\ast}(-\theta_{k}\theta_{j})=-\theta_{i}^{\ast}(\theta_{k}\theta_{j})$.
Moreover, observe that the following elements are in the same equivalence class $\theta_{i}\theta_{j}\theta_{k}$, $\theta_{j}\theta_{k}\theta_{i}$ and $\theta_{k}\theta_{i}\theta_{j}$. Hence, they have the same image under $\theta^{\ast}$.\\

To illustrate this, let's compute for example $\theta_{4}^{\ast}(\theta_{2}\theta_{4} \theta_{7})$. 

$\begin{array}{lcl} \theta_{4}^{\ast}(\theta_{2}\theta_{4} \theta_{7})&=&\theta_{4}^{\ast}(\theta_{2})\theta_{4} \theta_{7} + (-1)^{\mid\theta_{2}\mid} \theta_{2}\theta_{4}^{\ast}(\theta_{4} \theta_{7})\\\\&=& 0 - \theta_{2} [\theta_{4}^{\ast}(\theta_{4})\theta_{7}+(-1)^{\mid\theta_{4}\mid}\theta_{4} \theta_{4}^{\ast}(\theta_{7})] \;since\;\theta_{4}^{\ast}(\theta_{2})=0\\\\&=&
-\theta_{2} \theta_{7}\;since\;\theta_{4}^{\ast}(\theta_{7})=0,\mid\theta_{4}\mid\;=\;1\;and\;\theta_{4}^{\ast}(\theta_{4})=1

\end{array}$\\
%
%
\\
In order to complete the description of $d$, we need to say what $\partial$ is.
$$\partial_{i}: k[x_{1},\cdots,x_{n},x_{1}^{'},\cdots,x_{n}^{'}]\longrightarrow k[x_{1},\cdots,x_{n},x_{1}^{'},\cdots,x_{n}^{'}]$$ is defined by,

$$\partial_{i}(h)=\frac{h(x_{1}^{'},\cdots,x_{i-1}^{'},x_{i},\cdots,x_{n},x_{1}^{'},\cdots,x_{n}^{'})-h(x_{1}^{'},
\cdots,x_{i}^{'},x_{i+1},\cdots,x_{n},x_{1}^{'},\cdots,x_{n}^{'})}{x_{i}-x_{i}^{'}}$$
where for ease of notation we wrote $h$ as argument of $\partial_{i}$ instead of the more cumbersome notation
$\partial_{i}(h(x_{1},\cdots,x_{n},x_{1}^{'},\cdots,x_{n}^{'}))$, we will do same in what follows. But first, observe that $\partial_{i}$ is well defined because its numerator will always have $(x_{i}-x_{i}^{'})$ as a factor. Thus, ensuring that $\partial_{i}$ is a polynomial.

In the sequel, we will sometimes write $x$ for $x_{1},\cdots,x_{n}$ and $x^{'}$ for $x_{1}^{'},\cdots,x_{n}^{'}$.\\

Just like in \cite{carqueville2016adjunctions}, we can write: $$\partial_{i}(h)=\frac{^{t_{1}\cdots t_{i-1}}h(x_{1},\cdots,x_{n},x_{1}^{'},\cdots,x_{n}^{'}) - ^{t_{1}\cdots t_{i}}h(x_{1},\cdots,x_{n},x_{1}^{'},\cdots,x_{n}^{'})}{x_{i}-x_{i}^{'}}$$

where $^{t_{i}}(-): K[x,x'] \rightarrow K[x,x']$ , $h\mapsto h|_{x_{i}\mapsto x_{i}^{'}}$\\

is a variable changing map which in any polynomial, replaces the variable $x_{i}$ by the variable $x'_{i}$.\\
So, in particular for the $f$ in $\Delta_{f}$, which is also the same $f$ used in defining the differential map $d$ (see $\natural$ above), we have  $f\in K[x]\subseteq K[x,x']$ and $$\partial_{i}(f)=\frac{^{t_{1}\cdots t_{i-1}}f(x_{1},\cdots,x_{n}) - ^{t_{1}\cdots t_{i}}f(x_{1},\cdots,x_{n})}{x_{i}-x_{i}^{'}}$$
\begin{example}
Let $f=x-y \in \mathbb{R}[x,y]$.
  $$\partial_{1}(f(x,y))=\cfrac{f(x,y)-f(x',y)}{x-x'}=\cfrac{(x-y)-(x'-y)}{x-x'}=1.$$
   $$\partial_{2}(f(x,y))=\cfrac{f(x',y)-f(x',y')}{y-y'}=\cfrac{(x'-y)-(x'-y')}{y-y'}=-1.$$
\end{example}
We now prove the following lemma whose proof is omitted in \cite{carqueville2016adjunctions}.
\begin{lemma}\cite{carqueville2016adjunctions}
For\footnote{the $f$ in this lemma should not be confused with the $f$ before this lemma or after this lemma.} $f,g\in K[x,x']$, we have $$\partial_{i}(fg)= \partial_{i}(f)(^{t_{1}\cdots t_{i}} g) + (^{t_{1}\cdots t_{i-1}}f)\partial_{i}(g).$$

\end{lemma}
\begin{proof}
First of all, observe that by definition of $^{t_{i}}(-)$, it is obvious that $$^{t_{1}\cdots t_{i}}(fg) =(^{t_{1}\cdots t_{i}}f )(^{t_{1}\cdots t_{i}}g). $$
$\begin{array}{lcl}
(\partial_{i}f(x,x^{'}))(^{t_{1}\cdots t_{i}}g(x,x^{'}))+ (^{t_{1}\cdots t_{i-1}}f(x,x^{'}))(\partial_{i} g)\\\\=
\cfrac{^{t_{1}\cdots t_{i-1}}f(x,x^{'})^{t_{1}\cdots t_{i}}g(x,x^{'}) - ^{t_{1}\cdots t_{i}}f(x,x^{'})^{t_{1}\cdots t_{i}}g(x,x^{'})}{x_{i}-x_{i}^{'}} \\\\
+\cfrac{^{t_{1}\cdots t_{i-1}}f(x,x^{'})^{t_{1}\cdots t_{i-1}}g(x,x^{'}) - ^{t_{1}\cdots t_{i-1}}f(x,x^{'})^{t_{1}\cdots t_{i}}g(x,x^{'})}{x_{i}-x_{i}^{'}}
\\\\= \cfrac{^{t_{1}\cdots t_{i-1}}f(x,x^{'})^{t_{1}\cdots t_{i-1}}g(x,x^{'}) - ^{t_{1}\cdots t_{i}}f(x,x^{'})^{t_{1}\cdots t_{i}}g(x,x^{'})}{x_{i}-x_{i}^{'}}
\\\\= \cfrac{^{t_{1}\cdots t_{i-1}}(f(x,x^{'})g(x,x^{'})) - ^{t_{1}\cdots t_{i}}(f(x,x^{'})g(x,x^{'}))}{x_{i}-x_{i}^{'}}
\\\\=\partial_{i}(fg)\\ as\,desired.
\end{array}$\\\\
The first two equalities hold by definition of $\partial_{i}$ .
\end{proof}
It is worth mentioning at this point in time that the authors in \cite{carqueville2016adjunctions} state what the differential $d$ for $\Delta_{f}$ is, but they do not prove that $(\Delta_{f}, d)$ determines a matrix factorization of
 $  f\otimes id - id\otimes f $ (and this is important in order to see that $(\Delta_{f}, d)$ is the unit matrix factorization with respect to the tensor product of matrix factorizations).\\
We now have enough ingredients to produce such a proof. We state the following lemma which shows that $(\Delta_{f}, d)$ determines a finite rank matrix factorization of $  f\otimes id - id\otimes f $. The proof is a bit lengthy.

\begin{lemma}
  The $R\otimes R-$module $\Delta_{f}$ together with the differential $d$ defined above, determine a finite rank matrix factorization of $ f\otimes id - id\otimes f$ (which is equivalent to both $h(x,x')$ and $f(x)-f(x')$).
\end{lemma}
\begin{proof}
As already discussed under subsection \ref{subsec: unit 1-morphisms in LGk}, $\Delta_{f}$ is a finite rank $R\otimes R-$module. That is why we will conclude by the end of this proof that $(\Delta_{f}, d)$ is a finite rank matrix factorization given that its underlying module is of finite rank.\\
  Next, we need to show that $d$ is an odd degree map and that $d^{2}= h(x,x')\cdot id$, where $h(x,x')= f(x)-f(x')$\\
   $d$ is an odd degree map: In fact, since $\Delta_{f}$ is generated by
  the symbols $\theta_{1},\cdots,\theta_{n}$ it suffices to consider the action of $d$ on an arbitrary product of elements from the set $\{\theta_{1},\cdots,\theta_{n}\}$.\\
  \underline{\textit{Case 1}}:\\
   Let $p=\theta_{j_{1}}\theta_{j_{2}}\cdots \theta_{j_{k}}$ be an odd degree element, i.e.,
  $k$ is odd. We also have $\theta_{j_{l}}\in \{\theta_{1},\cdots,\theta_{n}\}$ with $l\in \{1,\cdots, k\}$. \\

  Claim 1: For each $i$, $(x_{i}-x_{i}')\theta_{i}^{\ast}(p)$ is either of even degree or is zero. This is easy to see as it is a direct consequence of the definition of $\theta_{i}^{\ast}$. In fact, we know by definition of $\theta_{i}^{\ast}$ that when it is applied to $p$, if $p$ contains no $\theta_{i}$, then the result would be zero. Now, if $p$ contains $\theta_{i}$, then this $\theta_{i}$ will no longer appear among the $\theta$s in $\theta_{i}^{\ast}(p)$, but the other $\theta$s that were in $p$ will remain. Thus, an even number of $\theta$s will remain in $\theta_{i}^{\ast}(p)$. For e.g. we saw after computation above that, $\theta_{4}^{\ast}(\theta_{2}\theta_{4}\theta_{7})=-\theta_{2}\theta_{7}$ hence $(x_{4}-x_{4}')\theta_{4}^{\ast}(\theta_{2}\theta_{4}\theta_{7})=-(x_{4}-x_{4}')\theta_{2}\theta_{7}$ which is of even degree. We also have for e.g. $\theta_{4}^{\ast}(\theta_{2})=0$ hence $(x_{4}-x_{4}')\theta_{4}^{\ast}(\theta_{2})=0$. \\

 Claim 2: For each $i$, $\partial_{i}(f)\theta_{i}\wedge(p)$ is of even degree.
  Observe that the degree of $\partial_{i}(f)\theta_{i}$ is deg$\theta_{i}=1$ since $\partial_{i}(f)$
contains no $\theta$. Thus, deg$[\partial_{i}(f)\theta_{i}\wedge(p)]=$deg$[\theta_{i}\wedge(p)]=$deg$(p)+1$.
So the claim is proved.\\
It now follows that when $p$ is of odd degree, $d(p)$ is the summation of even degree elements and so is of even degree. \\
\underline{\textit{Case 2}}:\\
Let $p=\theta_{j_{1}}\theta_{j_{2}}\cdots \theta_{j_{k}}$ be an even degree element, i.e.,
  $k$ is even. We also have $\theta_{j_{l}}\in \{\theta_{1},\cdots,\theta_{n}\}$ with $l\in \{1,\cdots, k\}$. \\

  Claim 3: For each $i$, $(x_{i}-x_{i}')\theta_{i}^{\ast}(p)$ is of odd degree or zero.\\
  A similar reasoning to that of claim 1 above can be applied to show that when $p$ is of even degree, $(x_{i}-x_{i}')\theta_{i}^{\ast}(p)$ is of odd degree or zero.\\

Claim 4: For each $i$, $\partial_{i}(h)\theta_{i}\wedge(p)$ is of odd degree.
  A similar reasoning to that of claim 2 can be applied to show that when $p$ is of even degree, $\partial_{i}(h)\theta_{i}\wedge(p)$ is of odd degree.\\

Hence, when $p$ is of even degree, $d(p)$ is the summation of odd degree elements and so is of odd degree.\\
Consequently, $d$ is an odd degree map.\\
The following definition will be needed in the sequel.
\begin{definition}
If $i\in \{1, \cdots, n\}$ and $\theta_{i}\,\in \,\theta_{1},\cdots,\theta_{n}$; we say that $\theta_{i}$ \textbf{is even} in $\theta_{1},\cdots,\theta_{n}$ if $\theta_{i}$ occupies an even position in $\theta_{1},\cdots,\theta_{n}$
\end{definition}
\begin{example}
$\theta_{5}$ is even in $\theta_{2}\theta_{5}\theta_{6}$ and $\theta_{6}$ is odd in $\theta_{2}\theta_{5}\theta_{6}$.
\end{example}

We now show that $d^{2}=h(x,x')\cdot id$ viz., $d^{2}=f(x)-f(x')\cdot id$. The author of \cite{khovanov2008matrix} simply writes $d^{2}= w$, omitting the "$\cdot id$". So, according to this notation, we want to show $d^{2}=h(x,x')=f(x)-f(x')$. \\
To that end, define $$A_{i}=(x_{i}-x_{i}')\theta_{i}^{\ast}\;and\;B_{i}=\partial_{i}(f)\theta_{i}\wedge(-)$$
$\begin{array}{lcl}
 d(d(H))\\=\sum_{i=1}^{n}(x_{i}-x_{i}')\theta_{i}^{\ast}(\sum_{j=1}^{n}(x_{j}-x_{j}')\theta_{j}^{\ast}(H)+ \partial_{j}(f)\theta_{j}\wedge (H))+\\
 \partial_{i}(f)\theta_{i}\wedge (\sum_{j=1}^{n}(x_{j}-x_{j}')\theta_{j}^{\ast}(H)+\partial_{j}(f)\theta_{j}\wedge (H))\\\\
 = (x_{1}-x_{1}')\theta_{1}^{\ast}((x_{1}-x_{1}')\theta_{1}^{\ast}(H)+\partial_{1}(f)\theta_{1}\wedge(H)+
 (x_{2}-x_{2}')\theta_{2}^{\ast}(H)+\\
 \partial_{2}(f)\theta_{2}\wedge(H)+\cdots +(x_{n}-x_{n}')\theta_{n}^{\ast}(H)+\partial_{n}(f)\theta_{n}\wedge(H))\\
 +\partial_{1}(f)\theta_{1}\wedge ((x_{1}-x_{1}')\theta_{1}^{\ast}(H)+\partial_{1}(f)\theta_{1}\wedge(H)+
 (x_{2}-x_{2}')\theta_{2}^{\ast}(H)+\\
 \partial_{2}(f)\theta_{2}\wedge(H)+\cdots +(x_{n}-x_{n}')\theta_{n}^{\ast}(H)+\partial_{n}(f)\theta_{n}\wedge(H))\\
 + (x_{2}-x_{2}')\theta_{2}^{\ast}((x_{1}-x_{1}')\theta_{1}^{\ast}(H)+\partial_{1}(f)\theta_{1}\wedge(H)+
 (x_{2}-x_{2}')\theta_{2}^{\ast}(H)+\\
 \partial_{2}(f)\theta_{2}\wedge(H)+\cdots +(x_{n}-x_{n}')\theta_{n}^{\ast}(H)+\partial_{n}(f)\theta_{n}\wedge(H))\\
 +\partial_{2}(f)\theta_{2}\wedge ((x_{1}-x_{1}')\theta_{1}^{\ast}(H)+\partial_{1}(f)\theta_{1}\wedge(H)+
 (x_{2}-x_{2}')\theta_{2}^{\ast}(H)+
 \end{array}$\\\\
  $\begin{array}{lcl}\partial_{2}(f)\theta_{2}\wedge(H)+\cdots +(x_{n}-x_{n}')\theta_{n}^{\ast}(H)+\partial_{n}(f)\theta_{n}\wedge(H))\\
 + \cdots +\\
 + (x_{n}-x_{n}')\theta_{n}^{\ast}((x_{1}-x_{1}')\theta_{1}^{\ast}(H)+\partial_{1}(f)\theta_{1}\wedge(H)+
 (x_{2}-x_{2}')\theta_{2}^{\ast}(H)+\\
 \partial_{2}(f)\theta_{2}\wedge(H)+\cdots +(x_{n}-x_{n}')\theta_{n}^{\ast}(H)+\partial_{n}(f)\theta_{n}\wedge(H))\\
 +\partial_{n}(f)\theta_{n}\wedge ((x_{1}-x_{1}')\theta_{1}^{\ast}(H)+\partial_{1}(f)\theta_{1}\wedge(H)+
 (x_{2}-x_{2}')\theta_{2}^{\ast}(H)+\\
 \partial_{2}(f)\theta_{2}\wedge(H)+\cdots +(x_{n}-x_{n}')\theta_{n}^{\ast}(H)+\partial_{n}(f)\theta_{n}\wedge(H))

 \end{array}$
\\\\
Recalling that $$A_{i}=(x_{i}-x_{i}')\theta_{i}^{\ast}\;and\;B_{i}=\partial_{i}(f)\theta_{i}\wedge(-)$$
we get:\\\\
$ d(d(-))\\=A_{1}\circ A_{1} + A_{1}\circ B_{1}+A_{1}\circ A_{2} +A_{1}\circ B_{2}+\cdots+A_{1}\circ A_{n}+A_{1}\circ B_{n}\\
+B_{1}\circ A_{1} + B_{1}\circ B_{1}+B_{1}\circ A_{2} +B_{1}\circ B_{2}+\cdots+B_{1}\circ A_{n}+B_{1}\circ B_{n}\\
+A_{2}\circ A_{1} + A_{2}\circ B_{1}+A_{2}\circ A_{2} +A_{2}\circ B_{2}+\cdots+A_{2}\circ A_{n}+A_{2}\circ B_{n}\\
 +B_{2}\circ A_{1} + B_{2}\circ B_{1}+B_{2}\circ A_{2} +B_{2}\circ B_{2}+\cdots+B_{2}\circ A_{n}+B_{2}\circ B_{n}\\
 +\cdots \\
 +A_{n}\circ A_{1} + A_{n}\circ B_{1}+A_{n}\circ A_{2} +A_{n}\circ B_{2}+\cdots+A_{n}\circ A_{n}+A_{n}\circ B_{n}\\
 +B_{n}\circ A_{1} + B_{n}\circ B_{1}+B_{n}\circ A_{2} +B_{n}\circ B_{2}+\cdots+B_{n}\circ A_{n}+B_{n}\circ B_{n}
$\\\\
Thus we get: \\

$$d^{2}=\sum_{i,j=1}^{n}[A_{i}\circ A_{j}+A_{i}\circ B_{j}+B_{i}\circ A_{j}+B_{i}\circ B_{j}]$$
Observations:
\begin{itemize}
  \item For all $i$, $A_{i}\circ A_{i}=B_{i}\circ B_{i}=0$.\\
  We first show that $A_{i}\circ A_{i}=0$, i.e., $A_{i}\circ A_{i}(H)=0$ for all $H\in \Delta_{f}$,
  where $H=h(x,x')\theta_{j_{1}}\theta_{j_{2}}\cdots \theta_{j_{k}}=h(x,x')p$, with $p=\theta_{j_{1}}\theta_{j_{2}}\cdots \theta_{j_{k}}$; and $\theta_{j_{1}},\theta_{j_{2}},\cdots, \theta_{j_{k}}\in \{\theta_{1},\theta_{2},\cdots, \theta_{n}\}$. We have $(A_{i}\circ A_{i})(H)=h(x,x')(A_{i}\circ A_{i})(p)$ since $h(x,x')\in K[x,x']$ and $\theta^{\ast}$ as a derivation, is ($K[x,x']-$)linear\footnote{$\theta^{\ast}$ is an endomorphism on $\Delta_{f}$ which is an $R\otimes R$-module} i.e., a $K[x,x']$-module.\\
 We always have $(A_{i}\circ A_{i})(p)=0$. In fact, if $p$ contains $\theta_{i}$, then $A_{i}(p)$ will not contain any $\theta_{i}$, and so $(A_{i}\circ A_{i})(p)=0$. Now, if $p$ does not contain $\theta_{i}$, then $A_{i}(p)=0$, hence $(A_{i}\circ A_{i})(p)=A_{i}(0)=0$.

   Next, to prove $B_{i}\circ B_{i}=0$, recall that $\theta_{i}\wedge \theta_{i}=\theta_{i}\theta_{i}=0$.

   $B_{i}\circ B_{i}(H)=\partial_{i}(f)\theta_{i}\wedge(\partial_{i}(f)\theta_{i}\wedge(H))=
   \partial_{i}(f)\partial_{i}(f)(\theta_{i}\wedge \theta_{i})\wedge H=0$.
  \item If $i\neq j$, then $A_{i}\circ B_{j}=-B_{j}\circ A_{i}$
  Since $\Delta_{f}$ is generated by the $\theta_{i}$s, it suffices to verify that this equality holds for an arbitrary $p=\theta_{j_{1}}\theta_{j_{2}}\cdots\theta_{j_{k}}\in \Delta_{f}$. It is good to keep in mind the assumption $i\neq j$.\\

  $\begin{array}{lcl}
     A_{i}\circ B_{j}(\theta_{j_{1}}\theta_{j_{2}}\cdots\theta_{j_{k}}) & = &  A_{i}( \partial_{j}(f)\theta_{j}\wedge(\theta_{j_{1}}\theta_{j_{2}}\cdots\theta_{j_{k}}))\\\\
    & = & (x_{i}-x_{i}')\theta_{i}^{\ast}(\partial_{j}(f)\theta_{j}\wedge(\theta_{j_{1}}\theta_{j_{2}}\cdots\theta_{j_{k}})) \\\\
     & = & (x_{i}-x_{i}')\partial_{j}(f)\theta_{i}^{\ast}(\theta_{j}\wedge(\theta_{j_{1}}\theta_{j_{2}}\cdots\theta_{j_{k}}))
   \\\\
   & = & \begin{cases} 0, \;i\notin \{j_{1},\cdots,j_{k}\} \\

   (x_{i}-x_{i}')\partial_{j}(f)(\theta_{j}(\theta_{j_{1}}\theta_{j_{2}}\cdots\hat{\theta_{i}}\cdots\theta_{j_{k}})),\;\theta_{i}\;is\;even\;in\;\theta_{j_{1}}\cdots\theta_{j_{k}} \\
    -(x_{i}-x_{i}')\partial_{j}(f)(\theta_{j}(\theta_{j_{1}}\theta_{j_{2}}\cdots\hat{\theta_{i}}\cdots\theta_{j_{k}})),\;\theta_{i}\;is\;odd\;in\;\theta_{j_{1}}\cdots\theta_{j_{k}} \end{cases}
   \end{array}$

   Next, we compute:\\

   $\begin{array}{lcl}
   B_{j}\circ A_{i}(\theta_{j_{1}}\theta_{j_{2}}\cdots\theta_{j_{k}})& = & \partial_{j}(f)\theta_{j}\wedge ((x_{i}-x_{i}')\theta_{i}^{\ast}(\theta_{j_{1}}\theta_{j_{2}}\cdots\theta_{j_{k}}))\\\\
   & = &
   (x_{i}-x_{i}')\partial_{j}(f)\theta_{j}\wedge (\theta_{i}^{\ast}(\theta_{j_{1}}\theta_{j_{2}}\cdots\theta_{j_{k}}))\\\\
   & = &\begin{cases} 0, \;i\notin \{j_{1},\cdots,j_{k}\} \\

   -(x_{i}-x_{i}')\partial_{j}(f)(\theta_{j}(\theta_{j_{1}}\theta_{j_{2}}\cdots\hat{\theta_{i}}\cdots\theta_{j_{k}})),\theta_{i}\;is\;even\;in\;\theta_{j_{1}}\cdots\theta_{j_{k}} \\
    (x_{i}-x_{i}')\partial_{j}(f)(\theta_{j}(\theta_{j_{1}}\theta_{j_{2}}\cdots\hat{\theta_{i}}\cdots\theta_{j_{k}})),\;\theta_{i}\;is\;odd\;in\;\theta_{j_{1}}\cdots\theta_{j_{k}} \end{cases}
   \end{array}$\\

  Hence, $A_{i}\circ B_{j}(\theta_{j_{1}}\theta_{j_{2}}\cdots\theta_{j_{k}}) = -B_{j}\circ A_{i}(\theta_{j_{1}}\theta_{j_{2}}\cdots\theta_{j_{k}})$. So, $A_{i}\circ B_{j}=-B_{j}\circ A_{i}$ as desired.
  \item If $i\neq j$, then $A_{i}\circ A_{j}=-A_{j}\circ A_{i}$.\\
  We consider $p$ as before and verify that $A_{i}\circ A_{j}(p)=-A_{j}\circ A_{i}(p)$.\\
  First of all,
   if either $i$ or $j$ is not in $I=\{j_{1},\cdots,j_{k}\}$,
  then $A_{i}\circ A_{j}=0=A_{j}\circ A_{i}$, and so $A_{i}\circ A_{j}=-A_{j}\circ A_{i}$ as desired.\\
  Next, suppose that $i,j\in I$ and $j\leq i$ without loss of generality. We distinguish four cases.\\
  \underline{\textit{Case 1}}: $\theta_{i}$ and $\theta_{j}$ are both even in $\{\theta_{j_{1}}\theta_{j_{2}}\cdots\theta_{j_{k}}\}$\\

  $\begin{array}{lcl}
   A_{i}\circ A_{j}(\theta_{j_{1}}\theta_{j_{2}}\cdots\theta_{j_{k}})  & = & (x_{i}-x_{i}')\theta_{i}^{\ast}((x_{j}-x_{j}')\theta_{j}^{\ast}(\theta_{j_{1}}\theta_{j_{2}}\cdots\theta_{j_{k}})) \\\\
     & = & -(x_{i}-x_{i}')(x_{j}-x_{j}')\theta_{i}^{\ast}(\theta_{j_{1}}\theta_{j_{2}}\cdots\hat{\theta_{j}}\cdots\theta_{i}\cdots\theta_{j_{k}}),\,\theta_{j}\,is\,even\,in\,\theta_{j_{1}}\theta_{j_{2}}\cdots\theta_{j_{k}}
\\\\
     & = & -(x_{i}-x_{i}')(x_{j}-x_{j}')(\theta_{j_{1}}\theta_{j_{2}}\cdots\hat{\theta_{j}}\cdots\hat{\theta_{i}}\cdots\theta_{j_{k}}),\,\theta_{i}\,is\,odd\,in\,\theta_{j_{1}}\cdots\hat{\theta_{j}}\cdots\theta_{i}\cdots\theta_{j_{k}}

  \end{array}$\\

  Next, we compute:\\

  $\begin{array}{lcl}
   A_{j}\circ A_{i}(\theta_{j_{1}}\theta_{j_{2}}\cdots\theta_{j_{k}}) & = & (x_{j}-x_{j}')\theta_{j}^{\ast}((x_{i}-x_{i}')\theta_{i}^{\ast}(\theta_{j_{1}}\theta_{j_{2}}\cdots\theta_{j_{k}})) \\\\
     & = & -(x_{j}-x_{j}')\theta_{j}^{\ast}(x_{i}-x_{i}')(\theta_{j_{1}}\cdots\theta_{j}\cdots\hat{\theta_{i}}\cdots\theta_{j_{k}}),\,\theta_{i}\,is\,even\,in\,\theta_{j_{1}}\theta_{j_{2}}\cdots\theta_{j_{k}}
\\\\
     & = & (x_{j}-x_{j}')(x_{i}-x_{i}')(\theta_{j_{1}}\cdots\hat{\theta_{j}}\cdots\hat{\theta_{i}}\cdots\theta_{j_{k}}),\,\theta_{j}\,is\,even\,in\,\theta_{j_{1}}\cdots\theta_{j}\cdots\hat{\theta_{i}}\cdots\theta_{j_{k}}

  \end{array}$\\

  It follows that
  $(A_{j}\circ A_{i})(\theta_{j_{1}}\theta_{j_{2}}\cdots\theta_{j_{k}})=-(A_{i}\circ A_{j})(\theta_{j_{1}}\theta_{j_{2}}\cdots\theta_{j_{k}})$, thus $A_{j}\circ A_{i}=-A_{i}\circ A_{j}$ as desired.\\

  \underline{Case \textit{2}}: $\theta_{i}$ and $\theta_{j}$ are both odd. \\
  The result holds here by following a reasoning completely similar to that of case 1.\\

 \underline{ \textit{Case 3}}: $\theta_{i}$ and $\theta_{j}$ are respectively odd and even.\\
  $\begin{array}{lcl}
   A_{i}\circ A_{j}(\theta_{j_{1}}\theta_{j_{2}}\cdots\theta_{j_{k}})  & = & (x_{i}-x_{i}')\theta_{i}^{\ast}((x_{j}-x_{j}')\theta_{j}^{\ast}(\theta_{j_{1}}\theta_{j_{2}}\cdots\theta_{j_{k}})) \\\\
     & = & -(x_{i}-x_{i}')(x_{j}-x_{j}')\theta_{i}^{\ast}(\theta_{j_{1}}\theta_{j_{2}}\cdots\hat{\theta_{j}}\cdots\theta_{i}\cdots\theta_{j_{k}}),\,\theta_{j}\,is\,even\,in\,\theta_{j_{1}}\theta_{j_{2}}\cdots\theta_{j_{k}}
\\\\
     & = & (x_{i}-x_{i}')(x_{j}-x_{j}')(\theta_{j_{1}}\theta_{j_{2}}\cdots\hat{\theta_{j}}\cdots\hat{\theta_{i}}\cdots\theta_{j_{k}}),\,\theta_{i}\,is\,even\,in\,\theta_{j_{1}}\cdots\hat{\theta_{j}}\cdots\theta_{i}\cdots\theta_{j_{k}}

  \end{array}$\\

  Next, we compute:\\

  $\begin{array}{lcl}
   A_{j}\circ A_{i}(\theta_{j_{1}}\theta_{j_{2}}\cdots\theta_{j_{k}}) & = & (x_{j}-x_{j}')\theta_{j}^{\ast}((x_{i}-x_{i}')\theta_{i}^{\ast}(\theta_{j_{1}}\theta_{j_{2}}\cdots\theta_{j_{k}})) \\\\
     & = & (x_{j}-x_{j}')\theta_{j}^{\ast}(x_{i}-x_{i}')(\theta_{j_{1}}\cdots\theta_{j}\cdots\hat{\theta_{i}}\cdots\theta_{j_{k}}),\,\theta_{i}\,is\,odd\,in\,\theta_{j_{1}}\theta_{j_{2}}\cdots\theta_{j_{k}}
\\\\
     & = & -(x_{j}-x_{j}')(x_{i}-x_{i}')(\theta_{j_{1}}\cdots\hat{\theta_{j}}\cdots\hat{\theta_{i}}\cdots\theta_{j_{k}}),\,\theta_{j}\,is\,even\,in\,\theta_{j_{1}}\cdots\theta_{j}\cdots\hat{\theta_{i}}\cdots\theta_{j_{k}}

  \end{array}$\\

 Thus, under this case too, $A_{j}\circ A_{i}=-A_{i}\circ A_{j}$ as desired.\\

  \underline{\textit{Case 4}}: $\theta_{i}$ and $\theta_{j}$ are respectively even and odd.\\
  This case is analogous to case 3.\\
  Hence, we always have $A_{j}\circ A_{i}=-A_{i}\circ A_{j}$.\\

  \item If $i\neq j$, then $B_{i}\circ B_{j}=-B_{j}\circ B_{i}$.\\
  In fact, since $\theta_{i}\theta_{j}=-\theta_{j}\theta_{i}$, we have:\\
  $B_{i}\circ B_{j}(-)= \partial_{i}(f)\theta_{i}\wedge (\partial_{j}(f)\theta_{j}\wedge (-))
  =-\partial_{i}(f)\partial_{j}(f)\theta_{j}\theta_{i}(-)=-B_{j}\circ B_{i}(-)$\\

  \item If $H=r_{1}\otimes r_{2}\theta_{i_{1}}\cdots\theta_{i_{n}}\in \Delta_{f}$, then either $A_{i}\circ B_{i}(H)=0$ or $B_{i}\circ A_{i}(H)=0$. In either case the result will be of the form $(x_{i}-x_{i}')\partial_{i}(f)H$.\\
  \underline{\textit{Case 1}}: $i\in \{i_{1},\cdots,i_{n}\}$  \\
   $A_{i}\circ B_{i}(H)=(x_{i}-x_{i}')\theta_{i}^{\ast}(\partial_{i}(f)\theta_{i}\wedge r_{1}\otimes r_{2}\theta_{i_{1}}\cdots\theta_{i_{n}})=0$. Since we will have the product $\theta_{i}\theta_{i}$ in this expression, and this product is clearly $0$ and $\theta_{i}^{\ast}(0)=0$ as $\theta_{i}^{\ast}$ is a derivation. So the result follows.\\
    Still under this case,\\

    $\begin{array}{lcl}
       B_{i}\circ A_{i}(r_{1}\otimes r_{2}\theta_{i_{1}}\cdots\theta_{i_{n}})& = &\partial_{i}(f)\theta_{i}\wedge ((x_{i}-x_{i}')\theta_{i}^{\ast}(r_{1}\otimes r_{2}\theta_{i_{1}}\cdots\theta_{i_{n}}))  \\\\
       & = & \begin{cases}
       -(x_{i}-x_{i}')\partial_{i}(f)\theta_{i}\wedge(r_{1}\otimes r_{2}\theta_{i_{1}}\cdots\hat{\theta_{i}}\cdots\theta_{i_{n}}),\theta_{i}\,is\,even\,in\,\theta_{i_{1}}\cdots\theta_{i_{n}} \\ (x_{i}-x_{i}')\partial_{i}(f)\theta_{i}\wedge(r_{1}\otimes r_{2}\theta_{i_{1}}\cdots\hat{\theta_{i}}\cdots\theta_{i_{n}}),\theta_{i}\,is\,odd\,in\,\theta_{i_{1}}\cdots\theta_{i_{n}} \end{cases}
       \\\\
        & = & \begin{cases}
       -(x_{i}-x_{i}')\partial_{i}(f)\wedge(-r_{1}\otimes r_{2}\theta_{i_{1}}\cdots\theta_{i}\cdots\theta_{i_{n}}) \\ (x_{i}-x_{i}')\partial_{i}(f)\wedge(r_{1}\otimes r_{2}\theta_{i_{1}}\cdots\theta_{i}\cdots\theta_{i_{n}}) \end{cases}
       \\\\
        &=&(x_{i}-x_{i}')\partial_{i}(f)H
     \end{array}$\\
     The sign in the first part (respectively second part) of the second to the last equality is due to the fact $\theta_{i}$ is moved an odd number of times (respectively an even number of times) to occupy its position, and each move is affected by a minus sign. \\

\underline{\textit{Case 2}}: $i\notin \{i_{1},\cdots,i_{n}\}$\\
     $B_{i}\circ A_{i}(r_{1}\otimes r_{2}\theta_{i_{1}}\cdots\theta_{i_{n}})=\partial_{i}(f)\theta_{i}\wedge (x_{i}-x_{i}')\theta_{i}^{\ast}(r_{1}\otimes r_{2}\theta_{i_{1}}\cdots\theta_{i_{n}})$. This yields $0$ since $\theta_{i}^{\ast}(r_{1}\otimes r_{2}\theta_{i_{1}}\cdots\theta_{i_{n}})=0$ as $i\notin \{i_{1},\cdots,i_{n}\}$.\\
     Still under this case, \\
     $\begin{array}{lcl}
     A_{i}\circ B_{i}(H)&=&(x_{i}-x_{i}')\theta_{i}^{\ast}(\partial_{i}(f)\theta_{i}\wedge r_{1}\otimes r_{2}\theta_{i_{1}}\cdots\theta_{i_{n}})\\\\
     &=&(x_{i}-x_{i}')\partial_{i}(f) r_{1}\otimes r_{2}\theta_{i}^{\ast}(\theta_{i}\theta_{i_{1}} \cdots\theta_{i_{n}})\\\\

     &=&(x_{i}-x_{i}')(\partial_{i}(f) r_{1}\otimes r_{2}\theta_{i_{1}}\cdots\theta_{i_{n}}),\,by\,definition\,of\,\theta_{i}^{\ast}\\\\
     &=&(x_{i}-x_{i}')\partial_{i}(f)H

     \end{array}$\\

    So, $A_{i}\circ B_{i}(H)=0$ or $B_{i}\circ A_{i}(H)=0$, and in either case the result is in the desired form, $(x_{i}-x_{i}')\partial_{i}(f)H$.
\end{itemize}

The foregoing work helps to simplify the expression of $d^{2}$ as follows:\\\\

$\begin{array}{lcl}
d^{2}= \sum_{i,j=1}^{n}[A_{i}\circ A_{j}+A_{i}\circ B_{j}+ B_{i}\circ A_{j}+ B_{i}\circ B_{j}]\\\\=
\sum_{j=1}^{n}A_{1}\circ A_{j}+A_{1}\circ B_{j}+ B_{1}\circ A_{j}+ B_{1}\circ B_{j}
 +A_{2}\circ A_{j}+A_{2}\circ B_{j}+ B_{2}\circ A_{j}+ B_{2}\circ B_{j}+\\
 \cdots +
 A_{n}\circ A_{j}+A_{n}\circ B_{j}+ B_{n}\circ A_{j}+ B_{n}\circ B_{j}\\\\=
 A_{1}\circ A_{1}+A_{1}\circ B_{1}+ B_{1}\circ A_{1}+ B_{1}\circ B_{1}
 +A_{1}\circ A_{2}+A_{1}\circ B_{2}+ B_{1}\circ A_{2}+ B_{1}\circ B_{2}+\\
 \cdots +
 A_{1}\circ A_{n}+A_{1}\circ B_{n}+ B_{1}\circ A_{n}+ B_{1}\circ B_{n}+

 A_{2}\circ A_{1}+A_{2}\circ B_{1}+ B_{2}\circ A_{1}+ B_{2}\circ B_{1}\\
 +A_{2}\circ A_{2}+A_{2}\circ B_{2}+ B_{2}\circ A_{2}+ B_{2}\circ B_{2}+
 \cdots +
 A_{2}\circ A_{n}+A_{2}\circ B_{n}+ B_{2}\circ A_{n}+ B_{2}\circ B_{n}+\\
 \cdots +
 A_{n}\circ A_{1}+A_{n}\circ B_{1}+ B_{n}\circ A_{1}+ B_{n}\circ B_{1}
 +A_{n}\circ A_{2}+A_{n}\circ B_{2}+ B_{n}\circ A_{2}+ B_{n}\circ B_{2}+\\
 \cdots +
 A_{n}\circ A_{n}+A_{n}\circ B_{n}+ B_{n}\circ A_{n}+ B_{n}\circ B_{n}\\\\=

A_{1}\circ B_{1}+ B_{1}\circ A_{1}+ A_{2}\circ B_{2}+ B_{2}\circ A_{2}+
\cdots+
A_{n}\circ B_{n}+ B_{n}\circ A_{n}\\\\=
\sum_{i=1}^{n}[A_{i}\circ B_{i}+ B_{i}\circ A_{i}]

\end{array}$\\\\

The second to the last equality above results from applying the identities obtained above namely
$A_{i}\circ A_{i}=B_{i}\circ B_{i}=0$, and if $i\neq j$ then $A_{i}\circ B_{j}=-B_{j}\circ A_{i}$, $A_{i}\circ A_{j}=-A_{j}\circ A_{i}$ and $B_{i}\circ B_{j}=-B_{j}\circ B_{i}$.

Thus the sum representation for $d^{2}$ reduces to:
$$d^{2}(H)=\sum_{i=1}^{n}[A_{i}\circ B_{i}(H) \;or\; B_{i}\circ A_{i}(H)]$$
Since $A_{i}\circ B_{i}=0 \;or\; B_{i}\circ A_{i}=0$.

 This sum is telescoping, in fact:\\

$\begin{array}{lcl}
 d^{2}=\sum_{i=1}^{n}[A_{i}\circ B_{i} \;or\; B_{i}\circ A_{i}]\\\\
 = \sum_{i=1}^{n}(x_{i}-x_{i}')\partial_{i}f \\\\
=(x_{1}-x_{1}')\partial_{1}f + (x_{2}-x_{2}')\partial_{2}f +\cdots + (x_{n-1}-x_{n-1}')\partial_{n-1}f + (x_{n}-x_{n}')\partial_{n}f \\\\
= f(x_{1},\cdots, x_{n})- f(x_{1}',x_{2},\cdots, x_{n})\\
+ f(x_{1}',x_{2},\cdots, x_{n}) - f(x_{1}',x_{2}',x_{3},\cdots, x_{n})\\
+\cdots\\
+f(x_{1}',\cdots, x_{n-2}',x_{n-1},x_{n})-
f(x_{1}',\cdots,x_{n-1}',x_{n})\\
+f(x_{1}',\cdots, x_{n-1}',x_{n}) - f(x_{1}',\cdots, x_{n}')\\\\
=
f(x_{1},\cdots, x_{n}) - f(x_{1}',\cdots, x_{n}')\\\\
=f(x)-f(x')

\end{array}$\\

So $d^{2}=(f(x)-f(x'))\cdot id$ as desired.

\end{proof}

We call $\Delta_{f}$ the unit matrix factorization as it is the unit with respect to the tensor
product of matrix factorizations. It is also referred to as the stabilised diagonal \cite{dyckerhoff2011compact}
or Koszul model of the diagonal \cite{carqueville2016adjunctions}. The diagonal here refers to $R$ as an $R\otimes R-$module, with multiplication giving the module structure, which is a matrix factorization of $f\otimes id - id\otimes f$
(or equivalently $f(x)-f(x')$) with differential zero. We give it the zero differential structure, which works, since $f(x)-f(x')$ acts as $0$ on this module. Thus, we have a canonical object in $hmf(R\otimes R, f(x)-f(x'))$ namely $(R,0)$. \\
Since we are constructing a bicategory, we would like to verify at this level that
for $f\in R=K[x]$ and for each object $(R,f)$ of $\mathcal{LG}_{K}$, there is a functor $I_{(R,f)}$ that returns the identity $\Delta_{f}$ on
$(R,f)$.\\
Define $I_{(R,f)}: 1=\{\ast\}\longrightarrow hmf(R\otimes R, id\otimes f- f\otimes id)^{\omega}$, $\ast \mapsto \Delta_{f}$.\\
 $I_{(R,f)}$ is a functor. In fact:\\
 Clearly $\Delta_{f}$ is an object of $hmf(R\otimes R, id\otimes f- f\otimes id)$ as discussed above.\\
 Next, let $1_{\ast}: \ast \longrightarrow \ast$, then clearly $I_{(R,f)}(1_{\ast})= 1_{\Delta_{f}}= 1_{I_{(R,f)}(\ast)}$.\\
 Now, consider two maps $F$ and $G$ in the small category $\{\ast\}$, which is the final object in the category Cat of small categories, then clearly the only possibility for the two maps is $F=G=1_{\ast}$. We clearly have:\\

 $I_{(R,f)}(F\circ G)=I_{(R,f)}(1_{\ast}\circ 1_{\ast})=I_{(R,f)}(1_{\ast})=1_{\Delta_{f}}$ ... (i)\\
 $I_{(R,f)}(F)\circ I_{(R,f)}(G)=I_{(R,f)}(1_{\ast})\circ I_{(R,f)}(1_{\ast})=1_{\Delta_{f}} \circ 1_{\Delta_{f}}= 1_{\Delta_{f}}$... (ii)

 (i) and (ii) show that $I_{(R,f)}(F\circ G)=I_{(R,f)}(F)\circ I_{(R,f)}(G)$, and this completes the proof that $I_{(R,f)}$ is a functor.\\

 There is an important observation made by \cite{carqueville2016adjunctions} that we state here as lemma \ref{projection map} below.

\begin{lemma} \label{projection map}

There is a canonical map of factorizations \\$\pi: \Delta_{f}\longrightarrow R$ given by
$\pi[(r\otimes r')\theta_{i_{1}}\theta_{i_{2}}\cdots \theta_{i_{k}}]=\delta_{k,0}rr'$.
$\pi$ is in fact the composition of the projection $\pi^{\ast}: \Delta_{f}\longrightarrow R\otimes R$ to $\theta-$degree $0$, followed by the multiplication map $m: R\otimes R \longrightarrow R=R_{0}\oplus R_{1}$, where we endow $R$ with the trivial grading i.e., $R_{0}=R$ and $R_{1}=\{0\}$.
\end{lemma}

\begin{proof}

To show that $\pi$ is a map of factorizations, we need to show that it is grade preserving, $R\otimes R-$linear and satisfies
$d_{R}\pi=\pi d$, where $d$ and $d_{R}$ are respectively the differentials given to $\Delta_{f}$ and $R$, and $d_{R}=0$ as explained above.

\begin{enumerate}
  \item $\pi$ is an even map, i.e., a grade preserving map.\\
  It suffices to show that $\pi$ sends elements of even (resp. odd) degree in $\Delta_{f}$ to elements of even (resp. odd) degree in $R$. \\
  In fact, since the degree of the empty word (i.e., the word $\theta_{i_{1}}\theta_{i_{2}}\cdots \theta_{i_{k}}$, where $k<1$) is zero, we have that $r\otimes r'$ is of degree $0$ and its image is also of degree $0$ because $\pi(r\otimes r')=rr'\in R=R_{0}$.\\
  Next, we see by definition of $\pi$ that all odd degree elements are mapped to $0\in R_{1}=\{0\}$.\\
  Moreover, if $k\neq 0$ and is even, then $(r\otimes r')\theta_{i_{1}}\theta_{i_{2}}\cdots \theta_{i_{k}}$ is of even degree and its image under $\pi$ is also of even degree since $\pi[(r\otimes r')\theta_{i_{1}}\theta_{i_{2}}\cdots \theta_{i_{k}}]=0\in R=R_{0}$
  \\
  So, $\pi$ is a degree preserving map.

  \item $\pi$ is $R\otimes R-$linear. Let $r_{1}\otimes r_{1}' \in R\otimes R$. \\
  We need to show that $\pi ((r_{1}\otimes r_{1}')(r\otimes r')\theta_{i_{1}}\cdots\theta_{i_{k}})= (r_{1}\otimes r_{1}')\pi((r\otimes r')\theta_{i_{1}}\cdots\theta_{i_{k}})$:\\
  Now,\\
  if  $k\geq 1$ then:\\
  $\pi ((r_{1}\otimes r_{1}')(r\otimes r')\theta_{i_{1}}\cdots\theta_{i_{k}})=\pi ((r_{1}r\otimes r_{1}'r')\theta_{i_{1}}\cdots\theta_{i_{k}})= 0$\;\;\;\;\;...\;\;\;\;\;$\ast$\\
  and if $k < 1$ we have:\\
  $\pi ((r_{1}\otimes r_{1}')(r\otimes r')\theta_{i_{1}}\cdots\theta_{i_{k}})=\pi ((r_{1}r\otimes r_{1}'r')\theta_{i_{1}}\cdots\theta_{i_{k}})= r_{1}rr_{1}'r'$\;\;\;\;\;...\;\;\;\;\;$\ast\ast$\\
  Next,\\
  if  $k\geq 1$ then:\\
  $(r_{1}\otimes r_{1}')\pi((r\otimes r')\theta_{i_{1}}\cdots\theta_{i_{k}})= (r_{1}\otimes r_{1}')(0)=0$\;\;\;\;\;...\;\;\;\;\;$\ast '$\\
   and for $k < 1$ we have:\\
   $(r_{1}\otimes r_{1}')\pi((r\otimes r')\theta_{i_{1}}\cdots\theta_{i_{k}})= (r_{1}\otimes r_{1}')rr'=r_{1} r_{1}'rr'=r_{1}r r_{1}'r'$\;\;\;\;\;...\;\;\;\;\; $\ast\ast '$\\ This last equality is by commutativity in $R$, and the second to the last equality is due to the fact that the $R\otimes R-$module structure on $R$ is given by multiplication.\\
   Since $\ast$ and $\ast '$ are the same, and $\ast\ast$ is same as $\ast\ast '$, it follows that $$\pi ((r_{1}\otimes r_{1}')(r\otimes r')\theta_{i_{1}}\cdots\theta_{i_{k}})= (r_{1}\otimes r_{1}')\pi((r\otimes r')\theta_{i_{1}}\cdots\theta_{i_{k}})$$ as desired.\\
   To end the proof that $\pi$ is $R\otimes R-$linear, we need to show that it is additive.
   Now, this is true by definition of $\pi$.

\item Finally, we show that $d_{R}\pi =\pi d$.\\
  Recall that $d_{R}=0$ and so $d_{R}\pi =0$. Thus, it suffices to show that
  $\pi d=0$.\\
  Let $H=(r\otimes r')\theta_{i_{1}}\cdots\theta_{i_{k}}\in \Delta_{f}$.\\
  $\pi d(H)= \pi (\sum_{i=1}^{n}[(x_{i}-x_{i}')\theta_{i}^{\ast}(r\otimes r')\theta_{i_{1}}\cdots\theta_{i_{k}} + \partial_{i}(h)\theta_{i}\wedge (r\otimes r')\theta_{i_{1}}\cdots\theta_{i_{k}}])\\
  =\pi (\textcolor{brown}{(x_{1}-x_{1}')(r\otimes r')\theta_{i_{1}}\cdots\hat{\theta_{1}}\cdots\theta_{i_{k}} +\cdots +(x_{n}-x_{n}')(r\otimes r')\theta_{i_{1}}\cdots\hat{\theta_{n}}\cdots\theta_{i_{k}}}\\
  + \textcolor{blue}{\partial_{1}(h)(r\otimes r')\theta_{1}\theta_{i_{1}}\cdots\theta_{i_{k}}+\cdots+ \partial_{n}(h)(r\otimes r')\theta_{n}\theta_{i_{1}}\cdots\theta_{i_{k}}})
  $
\end{enumerate}
The image under $\pi$ of the expression in blue is zero because for each summand:\\
Either $\theta_{i} \in \{\theta_{i_{1}}\cdots\theta_{i_{k}}\}$ and so we will have $\theta_{i}^{2}$ in that summand which will cause the summand to boil down to zero
or $\theta_{i} \notin \{\theta_{i_{1}}\cdots\theta_{i_{k}}\}$ and so the image under $\pi$ of what is left will be zero by definition of $\pi$ since $n\geq 1$.\\

The image under $\pi$ of the expression in brown is zero because for each summand:\\
If in the summand there is a $\theta_{i_{l}}$, $(l\in \{1,\cdots,k\})$, then by definition of $\pi$, the image would be zero.\\
If in the summand there is no $\theta_{i_{l}}$, $(l\in \{1,\cdots,k\})$, then its image under $\pi$ is
$\pi ((x_{i}-x_{i}')(r\otimes r'))=(r\otimes r')\pi(x_{i}-x_{i}')$ which is zero because
$\pi(x_{i}-x_{i}')=
\pi(x_{i})- \pi(x_{i}')=0$ in $R$.
\end{proof}
We now work out $(\Delta_{f}, d_{\Delta_{f}})$ for a specific $f$. We will simply write $d$ for $d_{\Delta_{f}}$.
\begin{example}
Let $f=x\in \mathbb{R}[x]$ be a one variable polynomial over the polynomial ring $\mathbb{R}[x]$, where $\mathbb{R}$ is the set of real numbers. Let $R=\mathbb{R}[x]$.
We have the $R\otimes R-$module:
 $$\Delta_{f}=\bigwedge\{\bigoplus(R\otimes R)\theta\}$$

 This is an exterior algebra generated by the variable $\theta$.

 To be more precise, it is the free $R\otimes R-$module generated by the empty list and the symbol $\theta$, modulo the relation that the $\theta^{2}=0$ that is $\theta\wedge \theta=0$.
 Here, the "$\wedge$" product is taken over $K$ just like the tensor product.
 A typical element of $\Delta_{f}^{1}$ is $(r\otimes r')\theta$ or equivalently $h(x,x')\theta$ where $h(x,x')$ $\in$
 $\mathbb{R}[x,x']$. A typical element of $\Delta_{f}^{0}$ is $(r\otimes r')$ or equivalently $h(x,x')$ where $h(x,x')$ $\in$
 $\mathbb{R}[x,x']$.\\
The action of $R\otimes R$ is the obvious one. In fact, for example, if $r_{1}\otimes r_{2}\in R\otimes R$, its left and right actions on an element $(r\otimes r')\theta$ of $\Delta_{f}$
  simply yield $(r_{1}\otimes r_{2})(r\otimes r')\theta$ since $R$ is commutative.

  $\Delta_{f}$ is endowed with the $\mathbb{Z}_{2}-$grading given by $\theta-$degree (where deg$\theta=1$). Thus $deg\theta^{2}=deg 0= 0$.
 \\
  Next, we define the differential as follows:\\
  $$d: \Delta_{f}\longrightarrow \Delta_{f}$$
  $$d(-)=[(x-x^{'})\theta^{\ast}(-)+\partial(f)\theta \wedge (-)]$$
  $$\partial(f(x))=\cfrac{f(x)-f(x')}{x-x'}=\cfrac{x-x'}{x-x'}=1.$$
  We now show that $d^{0}\circ d^{1}= (f(x)-f(x'))\cdot id_{\Delta^{1}}$ and $d^{1}\circ d^{0}= (f(x)-f(x'))\cdot id_{\Delta^{0}}$. But first, $d^{0}
  :\Delta_{f}^{0}\longrightarrow \Delta_{f}^{1}$ and $d^{0}(h(x,x'))= (x-x^{'})\theta^{\ast}(h(x,x'))+\partial(x)\theta \wedge (h(x,x'))=\theta h(x,x')$ since $\partial(x)=\partial f(x)=1$ and $\theta^{\ast}(h(x,x'))=0$ as $\theta^{\ast}$ is a derivation w.r.t $\theta$.\\
  Next, $d^{1}
  :\Delta_{f}^{1}\longrightarrow \Delta_{f}^{0}$ and $d^{1}(h(x,x')\theta)= (x-x^{'})\theta^{\ast}(h(x,x')\theta)+\partial(x)\theta \wedge (h(x,x')\theta)=(x-x^{'})h(x,x')$ since $\theta^{\ast}(\theta)=1$ and $\theta \wedge (h(x,x')\theta)=0$ as $\theta\theta=0$.\\
  We now compute
  \\
  $\begin{array}{ccc}d^{0}\circ d^{1}(h(x,x')\theta)
  &=& d^{0}((x-x^{'})(h(x,x'))\\
  &=& (x-x^{'})\theta^{\ast}(x-x^{'})(h(x,x'))+\partial(x)\theta \wedge (x-x^{'})(h(x,x'))\\
  &=&\theta(x-x^{'})(h(x,x'))\,since\,\theta^{\ast}(x-x^{'})(h(x,x'))=0\,and\,\partial(x)=1\\
  &=&(x-x^{'})(h(x,x'))\theta\\
  &=& (x-x^{'})\cdot id_{\Delta^{1}}[(h(x,x'))\theta]
  \end{array}$\\
  i.e., $d^{0}\circ d^{1}=(x-x^{'})\cdot id_{\Delta^{1}} = (f(x)-f(x'))\cdot id_{\Delta^{1}}$
  \\\\
  and\\\\
  $\begin{array}{ccc}d^{1}\circ d^{0}(h(x,x'))
  &=& d^{1}((h(x,x')\theta)\\
  &=& (x-x^{'})\theta^{\ast}(h(x,x')\theta)+\partial(x)\theta \wedge (h(x,x')\theta)\\
  &=&(x-x^{'})(h(x,x'))\,since\,\theta^{\ast}(\theta)=1\,and\,\theta^{2}=0\\
  &=&(x-x^{'})(h(x,x'))\\
  &=& (x-x^{'})\cdot id_{\Delta^{0}}[(h(x,x'))]
  \end{array}$\\

  i.e., $d^{1}\circ d^{0}=(x-x^{'})\cdot id_{\Delta^{0}} = (f(x)-f(x'))\cdot id_{\Delta^{0}}$
\end{example}
\begin{example}
Let $f=x-y\in \mathbb{R}[x,y]$ be a two-variable polynomial over the polynomial ring $\mathbb{R}[x,y]$. Let $R=\mathbb{R}[x,y]$.
  We have the $R\otimes R-$module:
 $$\Delta_{f}=\bigwedge\{\bigoplus_{i=1}^{2}(R\otimes R)\theta_{i}\}$$

 This is an exterior algebra generated by $2$ anti-commuting variables, $\theta_{1}$ and $\theta_{2}$.

 To be more precise, for $i,j\in \{1,2\}$; it is the free $R\otimes R-$module generated by the set containing the empty list and the symbols $\theta_{i}$, modulo the relations that the $\theta$'s anti-commute, that is $\theta_{i}\wedge \theta_{j}=-\theta_{j}\wedge \theta_{i}$. Typically, we will omit the wedge product and write for instance $\theta_{i}\wedge \theta_{j}$ simply as $\theta_{i}\theta_{j}$. Here, the "$\wedge$" product is taken over $\mathbb{R}$ just like the tensor product.
 A typical element of $\Delta_{f}^{0}$ is either of the form $(r\otimes r')$ (or equivalently $h(x,y,x',y')$) or $(r\otimes r')\theta_{i}\theta_{j}$ (or equivalently $h(x,y,x',y')\theta_{i}\theta_{j}$) where $h(x,y,x',y')\in \mathbb{R}[x,y,x',y']$. A typical element of $\Delta_{f}^{1}$ is $(r\otimes r')\theta_{i}$ or equivalently $h(x,y,x',y')\theta_{i}$ where $h(x,y,x',y')\in \mathbb{R}[x,y,x',y']$.\\
The action of $R\otimes R$ is the obvious one. That is, for example, if $r_{1}\otimes r_{2}\in R\otimes R$, its left and right actions on an element $(r\otimes r')\theta_{i}\theta_{j}$ of $\Delta_{f}$
  simply yield $(r_{1}\otimes r_{2})(r\otimes r')\theta_{i}\theta_{j}=(r_{1}r\otimes r_{2} r')\theta_{i}\theta_{j}$ since $R$ is commutative.

  $\Delta_{f}$ is endowed with the $\mathbb{Z}_{2}-$grading given by $\theta-$degree (where deg$\theta_{i}=1$ for each $i$). Thus $deg\theta_{i}^{2}=0\;and\;deg\theta_{i}\theta_{j}=0$.
 \\
  Next, we define the differential as follows:\\
  $$d: \Delta_{f}\longrightarrow \Delta_{f}$$
  $$d(-)=\sum_{i=1}^{2}[(x_{i}-x_{i}^{'})\theta_{i}^{\ast}(-)+\partial_{i}(f)\theta_{i}\wedge (-)]$$
   $$\partial_{1}(f(x,y))=\cfrac{f(x,y)-f(x',y)}{x-x'}=\cfrac{(x-y)-(x'-y)}{x-x'}=1.$$
   $$\partial_{2}(f(x,y))=\cfrac{f(x',y)-f(x',y')}{y-y'}=\cfrac{(x'-y)-(x'-y')}{y-y'}=-1.$$
  We now need to show that $d^{0}\circ d^{1}=[(x-y)-(x'-y')]\cdot id_{\Delta^{1}}$ and $d^{1}\circ d^{0}=[(x-y)-(x'-y')]\cdot id_{\Delta^{0}}$. We will not show all the details because it is computed in a manner similar to the calculations in the above example. We only specify $d^{0}$ and $d^{1}$.\\
   $d^{0}
  :\Delta_{f}^{0}\longrightarrow \Delta_{f}^{1}$ and $d^{0}(-)= (x-x^{'})(\theta_{1}^{\ast})^{0}(-)+\partial_{1}(f)\theta_{1} \wedge (-)+(y-y^{'})(\theta_{2}^{\ast})^{0}(-)+\partial_{2}(f)\theta_{2} \wedge (-)$
   i.e., $d^{0}(-)= (x-x^{'})(\theta_{1}^{\ast})^{0}(-)+\theta_{1} \wedge (-)+(y-y^{'})(\theta_{2}^{\ast})^{0}(-)-\theta_{2} \wedge (-)$
   since $\partial_{1}(f)=1$ and $\partial_{2}(f)=-1$. \\
  Next, $d^{1}
  :\Delta_{f}^{1}\longrightarrow \Delta_{f}^{0}$ and $d^{1}(-)=(x-x^{'})(\theta_{1}^{\ast})^{1}(-)+\partial_{1}(f)\theta_{1} \wedge (-)+(y-y^{'})(\theta_{2}^{\ast})^{1}(-)+\partial_{2}(f)\theta_{2} \wedge (-)$
   i.e., $d^{0}(-)= (x-x^{'})(\theta_{1}^{\ast})^{1}(-)+\theta_{1} \wedge (-)+(y-y^{'})(\theta_{2}^{\ast})^{1}(-)-\theta_{2} \wedge (-)$
   since $\partial_{1}(f)=1$ and $\partial_{2}(f)=-1$.

\end{example}
\subsection{The left and the right units of $\mathcal{LG}_{K}$}
In this subsection, we define the left and right identities for our bicategory $\mathcal{LG}_{K}$.
After proving their naturality with respect to $2-$morphisms, we prove that they do not have direct inverses thereby justifying the fact that their inverses in \cite{carqueville2016adjunctions} are found only up to homotopy.\\
We denote the right (respectively left) unit by $\rho$ (respectively $\lambda$). In the whole of this subsection, we will be dealing only with the left unit $\rho$ and we will omit the statements and proofs of the results for $\lambda$ because they are similar to the ones presented for $\rho$.\\
Consider a $1$-morphism $X\in$ $hmf(R\otimes S, 1_{R}\otimes g - f\otimes 1_{S})^{\omega}=hmf(R\otimes S, id\otimes g - f\otimes id)^{\omega}$.\\
Thus, $X$ is a matrix factorization of $id\otimes g - f\otimes id$ and is also an $R\otimes S-$module.\\
Let $1_{X}:X \longrightarrow X$ be the identity map and $\pi$ be the projection defined in lemma \ref{projection map}.
\begin{remark}
  Observe that any $S-$module $N$ can be considered as an $R-$module by letting $rn:=f(r)n$ where $f: R\longrightarrow S$ is a homomorphism of commutative rings.
Observe that the $R\otimes S-$module $X$ can be considered as an $R-$module via the following $K-$homomorphism of commutative unitary rings $f: R\longrightarrow R\otimes S$ defined
 by $f(r)=r\otimes 1_{S}$,
and hence one can also see $X$ as an $R\otimes R-$module by means of the following multiplicative map which is a $K-$homomorphism of commutative unitary rings $m: R\otimes R\longrightarrow R$ defined by $m(r\otimes r')=rr'$.
Moreover, observe that the $R\otimes R-$module $\Delta_{f}$ can be considered as an $R-$module via the following homomorphism of commutative unitary rings $f: R\longrightarrow R\otimes R$ defined
by $f(r)=r\otimes 1_{R}$.\\
\end{remark}
Thanks to this remark, it makes sense to form the following tensor product over $R$: $X\otimes_{R} R$ and
$X\otimes_{R} \Delta_{f}$ since $X$ and $\Delta_{f}$ can be viewed as $R-$modules. Consequently, we will simply write $X\otimes R$ and $X\otimes \Delta_{f}$ for ease of notation.\\
Similarly, since the $R\otimes S-$module $X$ and the $S\otimes S-$module $\Delta_{g}$ can be viewed as $S-$modules, we can form the module $\Delta_{g} \otimes_{S} X$ that we simply write as $\Delta_{g} \otimes X$.
\\
We also consider the map $u: X\otimes R\longrightarrow X$ defined by $u(x\otimes r)=xr$. This definition makes sense since $X$ can be viewed as an $R-$module. $u$ is an isomorphism (See example 1 page 363 of \cite{dummit2004abstract}).\\

Now, define $\rho_{X}: X\otimes \Delta_{f}\longrightarrow X$ by $\rho_{X}:=u\circ (1_{X}\otimes \pi)$ and $\lambda_{X}: \Delta_{g}\otimes X\longrightarrow X$ by $\lambda_{X}:=u\circ (\pi\otimes 1_{X})$.
\\
$\rho_{X}$ and $\lambda_{X}$ are clearly morphisms in $hmf(R\otimes S, id\otimes g - f\otimes id)^{\omega}$.
\\
 We show that $\rho$ and $\lambda$ are natural w.r.t. $2-$morphisms in the variable $X$ and we explain why there is no direct inverse for $\rho_{X}$ and $\lambda_{X}$, for each $X$. In \cite{carqueville2016adjunctions}, it is proved that the left and right unit maps are isomorphisms up to homotopy, but the reason why they are not direct isomorphisms is not given explicitly, so we give that reason at the end of this paper.
 \begin{lemma}
   $\rho$ and $\lambda$ are natural w.r.t. $2-$morphisms in the variable $X$.
 \end{lemma}
 \begin{proof} The proof for $\rho$ is given and the one for $\lambda$ is omitted because it is similar.\\
 In fact:
 \begin{enumerate}
   \item $\rho$ is natural in $X$.\\ In fact,
   we have the functor $$F:=(-)\otimes \Delta_{f}: hmf(R\otimes S, id\otimes g - f\otimes id)^{\omega}\longrightarrow hmf(R\otimes S, id\otimes g - f\otimes id)^{\omega},$$ $X \mapsto X\otimes \Delta_{f}$. \\
   $F$ is well defined because $X\otimes \Delta_{f}$ is an $R\otimes S-$module since $X$ is an $R\otimes S-$module and $\Delta_{f}$ is an $R\otimes R-$module. \\
   We also have the identity functor: $$G:=Id: hmf(R\otimes S, id\otimes g - f\otimes id)^{\omega}\longrightarrow hmf(R\otimes S, id\otimes g - f\otimes id)^{\omega},$$ $X \mapsto X$.

   For each object $X\in hmf(R\otimes S, id\otimes g - f\otimes id)^{\omega}$, $\rho_{X}: X\otimes \Delta_{f}\longrightarrow X$ is a map in $hmf(R\otimes S, id\otimes g - f\otimes id)^{\omega}$.\\
   Components $\rho_{X}$ of $\rho$ are such that for each $p: X\longrightarrow Y$, $\rho_{Y}\circ F(p)= G(p)\circ \rho_{X}$, viz. $\rho_{Y}\circ p\otimes id= p\circ \rho_{X}$ because the following diagram obviously commutes:

   $\xymatrix{\ar@{}[dr]|{\color{black}\circlearrowright}
   X\otimes \Delta_{f} \ar[d]_{\rho_{X}} \ar[r]^{p\otimes id} & Y\otimes \Delta_{f} \ar[d]^{\rho_{Y}} \\
   Id(X)=X \ar[r]_{G(p)=p} & Id(Y)=Y }$
  (we wrote $id$ for $id_{\Delta_{f}}$)
   In fact, let $r\otimes r'\theta_{i_{1}}\cdots\theta_{i_{k}}\in \Delta_{f}$, $i_{1},\cdots,i_{k} \in \{1,\cdots,n\}$. Recalling that $\rho_{X}=u\circ (1_{X}\otimes \pi)$, $u: \xymatrix@1{X\otimes R\ar[r]^<(.3){\cong}& X}$, $x\otimes r\mapsto xr$; and also that $\pi$ is the $\theta-$degree $0$ map, we have on the one hand: \\

   $\begin{array}{ccc}
    (i)\;\;\; \rho_{Y}\circ (p\otimes id)(x\otimes (r\otimes_{K} r')\theta_{i_{1}}\cdots\theta_{i_{k}})& = & \rho_{Y}(p(x)\otimes (r\otimes_{K} r')\theta_{i_{1}}\cdots\theta_{i_{k}})\\
     &=& u(1_{Y}\otimes \pi)(p(x)\otimes (r\otimes_{K} r')\theta_{i_{1}}\cdots\theta_{i_{k}} ) \\& = &
      u(p(x)\otimes \pi((r\otimes_{K} r')\theta_{i_{1}}\cdots\theta_{i_{k}}))\\
        &= &\begin{cases}
        0\;if\;k\geq 1, \\
        u(p(x)\otimes rr')=p(x)rr',\;otherwise.
        \end{cases}
    \end{array}
   $\\

   On the other hand, we have:\\

   $\begin{array}{ccc}
   (ii)\;\;\; p\circ \rho_{X}(x\otimes (r\otimes_{K} r')\theta_{i_{1}}\cdots\theta_{i_{k}}) &=& p\circ u(1_{X}\otimes \pi)(x\otimes (r\otimes_{K} r')\theta_{i_{1}}\cdots\theta_{i_{k}})\\
   &=& p\circ u(x\otimes \pi((r\otimes_{K} r')\theta_{i_{1}}\cdots\theta_{i_{k}}))\\
   &=& \begin{cases} 0\;if\;k\geq 1, \\
   p(u(x\otimes rr'))=p(xrr')=p(x)rr'\;otherwise.
    \end{cases}
   \end{array}
   $\\

   The very last equality above is obtained by $R-$linearity of $p$ which is a map of factorization.\\
   The work done in (i) and (ii) shows that $\rho_{Y}\circ (p\otimes id)= p\circ \rho_{X}$. So $\rho$ is indeed a natural transformation.

   \item $\rho_{X}$ is a morphism of factorizations, viz. an $R-$linear even map commuting with the differentials.\\
       $\bullet$$\rho_{X}=u\circ (1_{X}\otimes \pi)$ is an $R-$linear even map.\\
       In fact, The projection map $\pi$ was shown to be linear and even, in the proof of lemma \ref{projection map}, and
       we know that the identity map $1_{X}$ is clearly linear and even, so $1_{X}\otimes \pi$ is a linear even map.\\
       Moreover, $u$ is a linear map (see example 1, page 363 of \cite{dummit2004abstract}). Hence $\rho_{X}$ is linear as composition of linear maps.\\
       It now remains to show that $u$ is an even map. We give the trivial grading to $R$, i.e., $R=R_{0}\oplus \{0\}$. Now, $(X_{0}\oplus X_{1})\otimes R = X_{0}\otimes R \oplus X_{1}\otimes R$. We have:\\
       $u: (X_{0}\oplus X_{1})\otimes R = X_{0}\otimes R \oplus X_{1}\otimes R\longrightarrow X_{0}\oplus X_{1}$ defined by: $u(x_{0}\otimes r_{0})=x_{0}r_{0}\in X_{0}$ and $u(x_{1}\otimes r_{0})=x_{1}r_{0}\in X_{1}$, since $X$ is an $R$-module. This shows that $u$ is an even map.  \\

       $\bullet$ $\rho_{X}$ commutes with the differentials.\\
       Let $Z= X\otimes \Delta_{f}$ and let's write each factorization with its pair of differentials as follows:\\
        $(X,d^{0}_{X},d^{1}_{X})$, $(Z,d^{0}_{Z},d^{1}_{Z})$ and $(\Delta_{f},d^{0}_{\Delta},d^{1}_{\Delta})$. Since $X$ and $Z=X\otimes \Delta_{f}$ are factorizations of $h:= id\otimes g - f\otimes id$, the following hold:

$$\begin{cases}
d^{1}_{X}d^{0}_{X}=h\cdot id_{X^{0}}\\
d^{0}_{X}d^{1}_{X}=h\cdot id_{X^{1}}\\
d^{1}_{Z}d^{0}_{Z}=h\cdot id_{Z^{0}}\\
d^{0}_{Z}d^{1}_{Z}=h\cdot id_{Z^{1}}
\end{cases}
$$
         We have the following pair of maps

         $\xymatrix {Z^{0}\ar [r]^{\rho^{0}} &X^{0}}$ and $\xymatrix {Z^{1}\ar [r]^{\rho^{1}} &X^{1}}$.\\
         To show that $\rho_{X}$ commutes with the differentials, it suffices to show that
          the following diagram commutes:\\

 $\xymatrix@ R=0.6in @ C=.75in{Z^{0} \ar[r]^{{d^{0}_{Z}}} \ar[d]_{\rho^{0}} &
Z^{1} \ar[d]^{\rho^{1}} \ar[r]^{{d^{1}_{Z}}} & Z^{0}\ar[d]^{\rho^{0}}\\
X^{0} \ar[r]^{{d^{0}_{X}}} & X^{1}\ar[r]^{{d^{1}_{X}}} & X^{0}}$\\

         i.e., $$
\begin{cases}
\rho^{0}d^{1}_{Z}= d^{1}_{X}\rho^{1}\;...\;(i)\\
d^{0}_{X}\rho^{0}= \rho^{1}d^{0}_{Z}\;...\;(ii)
\end{cases}
$$
It suffices
to prove only one of these equalities because they are equivalent\footnote{ We proved something quite similar using matrix representations in Remark \ref{eqce of equations of facto}.}.
\\
 We briefly discuss $(ii)$ without giving a full proof because of the amount of the details involved: \\
 First recall that $Z^{0}=(X\otimes \Delta_{f})^{0}=(X^{0}\otimes \Delta_{f}^{0})\oplus (X^{1} \otimes \Delta_{f}^{1}) $. \\
 Also recall that the $\mathbb{Z}_{2}-$graded $R\otimes S-$module $X\otimes \Delta_{f}$
is a finite rank matrix factorization of the polynomial $g-f$ and it comes with an odd $R\otimes S-$linear endomorphism $d_{Z}:Z \rightarrow Z$ s.t. $d_{Z}^{2}=(g-f)\cdot id_{Z}$. We have $$d_{Z}^{0}: Z^{0}=(X^{0}\otimes \Delta_{f}^{0})\oplus (X^{1}\otimes \Delta_{f}^{1})\rightarrow (X^{0}\otimes \Delta_{f}^{1})\oplus (X^{1}\otimes \Delta_{f}^{0})= Z^{1}.$$\\ Now, if for ease of notation we write $id_{\Delta}$ for $id_{\Delta_{f}}$ then $d_{Z}^{0}=(d_{X}\otimes id_{\Delta} +id_{X}\otimes d_{\Delta})^{0}$ where $id_{X}\otimes d_{\Delta}$ has the usual Koszul signs when applied to elements.
We have: $\rho^{0}=(u\circ (1_{X}\otimes \pi))^{0} =u^{0}\circ (1_{X}\otimes \pi)^{0} + u^{1}\circ (1_{X}\otimes \pi)^{1}$.\\
 $\rho^{1}=(u\circ (1_{X}\otimes \pi))^{1} =u^{1}\circ (1_{X}\otimes \pi)^{0} + u^{0}\circ (1_{X}\otimes \pi)^{1}$.\\
  $(1_{X}\otimes \pi)^{0}=(1_{X}^{0}\otimes \pi^{0}) \oplus (1_{X}^{1}\otimes \pi^{1})$, $(1_{X}\otimes \pi)^{1}=(1_{X}^{0}\otimes \pi^{1}) \oplus (1_{X}^{1}\otimes \pi^{0})$. \\
With all these data, the proof of ii) follows.
 \end{enumerate}
\end{proof}
The authors of \cite{carqueville2016adjunctions} give a difficult proof in which they exhibit the homotopy inverses of $\rho_{X}$ and $\lambda_{X}$ but they do not show why the direct inverses of $\rho_{X}$ and $\lambda_{X}$ do not exist.\\
In the following theorem, we show that $\rho_{X}$ and $\lambda_{X}$ do not have direct inverses. This helps to understand why $\rho_{X}$ and $\lambda_{X}$ can only have inverses up to homotopy.

\begin{theorem} \label{rho does not have a direct inverse}
 $\rho_{X}$ and $\lambda_{X}$ do not have direct inverses.
\end{theorem}
\begin{proof}
The proof for $\rho_{X}$ is presented and the one for $\lambda_{X}$ is omitted because it is similar.\\
We find a direct right inverse of $\rho_{X}$ that we call $\psi_{X}$, next we prove that $\psi_{X}$ is not a direct left inverse of $\rho_{X}$ and then proceed by contradiction to show that $\rho_{X}$ has no direct inverse.
 \\
$\bullet$ Define $\psi_{X}: X\longrightarrow X\otimes \Delta_{f}$ by $x\mapsto x\otimes 1_{R\otimes R}$.\\
First, it is important to show that $\psi$ is natural in $X$. It is easy to see that $\psi_{X}: X\longrightarrow X\otimes \Delta_{f}$ is a map in $hmf(R\otimes S, id\otimes g-f\otimes id)^{\omega}$. Let $X\,and\,Y\in hmf(R\otimes S, id\otimes g-f\otimes id)^{\omega}$ and $j:X\mapsto Y$. We need to show that the following diagram commutes:\\
 $\xymatrix{\ar@{}[dr]|{\color{black}\circlearrowright}
   Id(X)=X \ar[d]_{\psi_{X}} \ar[r]^{j} & Id(Y)=Y \ar[d]^{\psi_{Y}} \\
X\otimes \Delta_{f} \ar[r]_{j\otimes id} &Y\otimes \Delta_{f} }$

i.e., $\psi_{Y}\circ j= j\otimes id \circ \psi_{X}$.\\
Now, $\psi_{Y}\circ j(x)= j(x)\otimes 1_{R\otimes R}$ by definition of $\psi$. $\cdots$ $\dag$ \\
$j\otimes id \circ \psi_{X}(x)=j\otimes id(x\otimes 1_{R\otimes R})=j(x)\otimes id(1_{R\otimes R})=j(x)\otimes 1_{R\otimes R}$ $\cdots$ $\ddag$.\\
$\dag$ and $\ddag$ yield the desired equality. Hence $\psi$ is a natural transformation.\\
Now, $\psi_{X}=\psi_{X}^{0}\oplus \psi_{X}^{1}$, where for $x=x_{0}\oplus x_{1}$, we have:\\
$\psi_{X}^{0}: X^{0}\longrightarrow (X\otimes \Delta_{f})^{0}=(X^{0}\otimes \Delta_{f}^{0})\oplus (X^{1}\otimes \Delta_{f}^{1}))$, defined by $\psi_{X}^{0}(x_{0})= x_{0}\otimes 1_{R\otimes R}=x_{0}\otimes (1_{R}\otimes 1_{R})$\\
$\psi_{X}^{1}: X^{1}\longrightarrow (X\otimes \Delta_{f})^{1}=(X^{0}\otimes \Delta_{f}^{1})\oplus (X^{1}\otimes \Delta_{f}^{0})$, defined by $\psi_{X}^{1}(x_{1})=(x_{1}\otimes 1_{R\otimes R})=x_{1}\otimes (1_{R}\otimes 1_{R})$\\
$\bullet$  It is easy to see that $\psi_{X}$ is an $R-$linear map. $\psi_{X}$ is an even map by construction.\\
$\bullet$ We now show that $\psi_{X}$ commutes with the differentials. This condition is represented
diagrammatically by the commutativity of the following diagram where $Z= X\otimes \Delta_{f}$.\\

 $\xymatrix@ R=0.6in @ C=.75in{X^{0} \ar[r]^{{d^{0}_{X}}} \ar[d]_{\psi^{0}} &
X^{1} \ar[d]^{\psi^{1}} \ar[r]^{{d^{1}_{X}}} & X^{0}\ar[d]^{\psi^{0}}\\
Z^{0} \ar[r]^{{d^{0}_{Z}}} & Z^{1}\ar[r]^{{d^{1}_{Z}}} & Z^{0}}$\\

That is, we need to show that:

$$
\begin{cases}
\psi^{1}d^{0}_{X}=d^{0}_{Z} \psi^{0}\;...\;(i)' \\
\psi^{0}d^{1}_{X}=d^{1}_{Z}\psi^{1}\;...\;(ii)'
\end{cases}
$$
 It suffices to prove only one of these equalities because they are equivalent\footnote{ We proved something quite similar using matrix representations in Remark \ref{eqce of equations of facto}.}. \\

We prove $(i)'$: \\\\
$\begin{array}{ccc}
 d^{0}_{Z} \psi^{0}(x_{0})&=& (d_{X}^{0}\otimes id +id\otimes d_{\Delta}^{0})(x_{0}\otimes 1_{R\otimes R}) \\
    & = & d_{X}^{0}(x_{0})\otimes 1_{R\otimes R} + (-1)^{\mid x_{0} \mid} id(x_{0})\otimes d_{\Delta}^{0}(1_{R\otimes R}) \\
   & = & d_{X}^{0}(x_{0})\otimes 1_{R\otimes R} + x_{0}\otimes 0\;\;since\;d_{\Delta}^{0}(1_{R\otimes R})=0\;and\; \mid x_{0} \mid =0  \\
   &=& \psi^{1}d_{X}^{0}(x_{0})
 \end{array}$\\

 Hence $(i)'$ holds. So $\psi_{X}$ is a map of factorization.\\
$\bullet$ $\psi_{X}$ is a direct right inverse of $\rho_{X}$.\\
We need to show that $\rho_{X}\circ \psi_{X}=id_{X}$.\\
$\begin{array}{lcl}
   \rho_{X}\circ \psi_{X}(x_{0}\oplus x_{1}) & = & \rho_{X}(\psi_{X}(x_{0}) \oplus  \psi_{X}(x_{1})),\;\psi_{X}\;is\;linear \\
    &= & u(1_{X}\otimes \pi)(x_{0}\otimes 1_{R\otimes R} \oplus x_{1}\otimes 1_{R\otimes R})\;by\;definition\;of\;\rho_{X}\;and\;\psi_{X} \\
    & = & u((1_{X}\otimes \pi)(x_{0}\otimes 1_{R\otimes R}) \oplus (1_{X}\otimes \pi)(x_{1}\otimes 1_{R\otimes R}))\\
    &=& u(x_{0}\otimes 1_{R}1_{R} \oplus x_{1}\otimes 1_{R}1_{R})\;as\;\pi(1_{R\otimes R})=\pi(1_{R}\otimes 1_{_{R}})=1_{R}1_{R}.\\

   & = &  u(x_{0}\otimes 1_{R} \oplus x_{1}\otimes 1_{R})\\
   &=& u(x_{0}\otimes 1_{R}) \oplus u(x_{1}\otimes 1_{R}),\;since\;u\;is\;linear\\
   &= & x_{0}1_{R}\oplus x_{1}1_{R}\;by\;definition\;of\;u \\
   &=&  x_{0}\oplus x_{1}\;as\;desired
 \end{array}$\\

 So $\rho_{X}\circ \psi_{X}= id_{X}$.\\
i.e., $\psi_{X}$ is a direct right inverse of $\rho_{X}$. \\
$\bullet$ We now show that $\psi_{X}$ is not a (direct) left inverse to $\rho_{X}$. All we need show is that $\psi_{X}\circ \rho_{X}\neq id_{Z}$. Let $x\otimes (r\otimes r')\theta_{1} \in Z = X\otimes \Delta_{f}$,
then:\\

 $\begin{array}{lcl}
 \psi_{X}\circ \rho_{X}(x\otimes (r\otimes r')\theta_{1})
 &=&\psi_{X}(u(1_{X}\otimes \pi)(x\otimes (r\otimes r')\theta_{1}))\\
 &=& \psi_{X}(u(x\otimes \pi(r\otimes r')\theta_{1}))\\
 &=& \psi_{X}(0),\;since\;\pi((r\otimes r')\theta_{1})=0\\
 &=& 0\\
 &\neq& id_{Z}(x\otimes (r\otimes r')\theta_{1})

\end{array}$\\
$\bullet$ Finally, suppose towards a contradiction that $\rho_{X}$ has an inverse (i.e., a map that is both a right and a left direct iverse of $\rho_{X}$),
call it $\chi_{X}: X \mapsto X\otimes \Delta_{f}$.\\
Then: $\rho_{X}\circ \chi_{X}=id_{X}$ $\cdots$ $\sharp$ \\
and  $\chi_{X} \circ \rho_{X}=id_{X\otimes \Delta_{f}}$  $\cdots$ $\sharp'$  by definition of the inverse of a map.\\
And $\rho_{X} \circ \psi_{X}=id_{X}$  $\cdots$ $\sharp"$ since $\psi_{X}$ is a right inverse of $\rho_{X}$.\\
from $\sharp\,and\,\sharp"$ we have:\\
$\begin{array}{lcl}
 \rho_{X} \circ \psi_{X}=\rho_{X} \circ \chi_{X}
 &\Rightarrow& \chi_{X}\circ(\rho_{X} \circ \psi_{X})=\chi_{X}\circ(\rho_{X} \circ \chi_{X})\\
 &\Rightarrow& (\chi_{X}\circ\rho_{X}) \circ \psi_{X}=(\chi_{X}\circ\rho_{X}) \circ \chi_{X},\,by\,associativity\\
 &\Rightarrow& (id_{X\otimes \Delta_{f}}) \circ \psi_{X}=(id_{X\otimes \Delta_{f}}) \circ \chi_{X}\,since\,\chi\circ\rho_{X} = id_{X\otimes \Delta_{f}}\\
 &\Rightarrow& \psi_{X}=\chi_{X}
 \end{array}$\\
 But this last equality is a contradiction since $\psi_{X}$ which is not a left inverse of $\rho_{X}$ cannot be equal to $\chi_{X}$ which is a left inverse of $\rho_{X}$.\\
So $\rho_{X}$ is not invertible as claimed.
\end{proof}
 A similar work with the left unitor (or unit action) $\lambda_{X}: \Delta_{g}\otimes X \rightarrow X$ shows that it is not directly invertible.\\
 We refer the reader to sections 3 and 4 of \cite{carqueville2016adjunctions} for a proof showing that the unitors have inverses up to homotopy. The authors in \cite{carqueville2016adjunctions} used \textit{Atiyah classes} (cf. section 3 of \cite{carqueville2016adjunctions} ) in their proof.
 But they never explained why it does not work at a non-homotopic setting. We did it in the foregoing lemma.
\begin{remark}
 Observe that $\psi_{X}$ is a right homotopy inverse of $\rho_{X}$. In fact,
we need to verify that $\rho_{X}\circ \psi_{X}\sim id_{X}$ where $\sim$ stands for the homotopy relation.\\
(Recall that: $\rho_{X}\circ \psi_{X},  id_{X}: X\longrightarrow X$).\\
We need to find an odd degree $R$-linear map $\alpha: X\longrightarrow X$ s.t., $d_{X}\alpha + \alpha d_{X}=\rho_{X}\circ \psi_{X}-id_{X}$. It suffices to take $\alpha$ to be the zero map (which is clearly an $R$-linear odd degree map).\\
In fact, if $\alpha=0$ then:\\
$d_{X}\alpha + \alpha d_{X}=\rho_{X}\circ \psi_{X}-id_{X}\Longleftrightarrow 0=\rho_{X}\circ \psi_{X}-id_{X}$, i.e., $\rho_{X}\circ \psi_{X}=id_{X}$ which is true since $\psi_{X}$ is the right direct inverse of $\rho$.\\
A similar work can be done for the left unit map $\lambda_{X}$.
\end{remark}


\bibliography{fomatati_ref}
\addcontentsline{toc}{section}
{References}
\end{document}